\documentclass[a4paper,10pt]{amsart}
\usepackage{tikz}
\usetikzlibrary{matrix,arrows,decorations.pathmorphing}
\usepackage{tikz-cd}
\usepackage{amsfonts, amssymb, amsthm, amsmath}  
\usepackage{hyperref}
\usepackage{graphicx}
\usepackage{psfrag}
\newtheorem{thm}{Theorem}

\newtheorem{lemma}[thm]{Lemma}  
\newtheorem{remark}[thm]{Remark}  
\newtheorem{defn}[thm]{Definition}  
  
\newtheorem{claim}[thm]{Claim}  
  
\numberwithin{thm}{section}

\providecommand{\totl}[1]{\ensuremath{\lceil #1\rceil }}
\providecommand{\totb}[1]{\ensuremath{\underline{ #1}}}

\newcommand{\ex}{\mathbf}
\providecommand {\e}[1]{\mathfrak t^{#1}}

\providecommand{\C}[2]{\ensuremath {C^{#1,\underline{#2}}}}

\newcommand{\Msw}{\mathcal M^{st}}

\DeclareMathOperator{\id}{id}

\newcommand{\dbar}{\bar{\partial}}

\providecommand{\et}[2]{\ensuremath{\bold T^{#1}_{#2}}}

\providecommand{\abs}[1]{\left\lvert #1\right\rvert}

\author{Brett Parker   }
  
\title{Universal tropical structures for curves in exploded manifolds}  

\thanks{This work was supported by ARC grant DP1093094. AI models Gemini and ChatGPT were used to assist editing the second version of the paper to improve exposition.}

\begin{document}

\begin{abstract} For any stable curve $f$ in an exploded manifold, this paper constructs a family of curves $\hat f$ with universal tropical structure which contains $f$. Such a family has the property that any other family of curves containing $f$ is locally a small modification of a family which factors through $\hat f$. As such, families of curves with universal tropical structure play an important role in the analysis of the moduli stack of curves and the construction of Gromov--Witten invariants on exploded manifolds.

\end{abstract}

\maketitle

\tableofcontents

\section{Introduction}

This paper is one of a series of papers defining Gromov--Witten invariants for
exploded manifolds~\cite{iec,scgp}, a framework that ultimately yields new
tropical gluing theorems for Gromov--Witten invariants of symplectic manifolds.
Exploded manifolds combine small-scale smooth geometry with piecewise-linear,
or tropical, large-scale geometry.

Once the tropical data is fixed, exploded geometry behaves much like ordinary
differential geometry.  The new feature is that this tropical data is itself
part of the geometry.  It affects, for example, the deformation theory of
holomorphic curves, and even basic extension questions for maps defined on
subsets of exploded manifolds.  Moreover, this large-scale tropical geometry is
essential for gluing theorems arising from normal-crossings degenerations of
symplectic manifolds.

\smallskip
The purpose of this paper is to construct the universal tropical data governing
deformations of curves in exploded manifolds.  The first step is purely
tropical.  Associated to a curve \(f\) in an exploded manifold is a tropical
curve \(\totb f\), obtained from the tropical structure of \(f\).  Deformations
of \(\totb f\) are encoded by extensions of this tropical structure.  We prove
that there is a universal such extension: every other extension is obtained
uniquely as its pullback along an integral-affine map; see Theorem~\ref{mts}.
Concretely, this universal extension is cut out by integral-affine consistency
conditions on the possible positions of the vertices and lengths of the
internal edges of \(\totb f\).  When the target is a smooth manifold, this
reduces to the usual orthant of smoothing parameters for the nodes.

We then realize this universal tropical extension geometrically.  For any
stable curve \(f\) in an exploded manifold, Lemma~\ref{uts family} constructs a
family of curves \(\hat f\) whose tropical structure restricts at \(f\) to the
universal extension.  The main existence theorem, Theorem~\ref{G uts family},
gives such a family with the automorphism properties needed for the local
theory.  Lemma~\ref{uts map} then gives the corresponding local modelling property:
any other family \(\hat h\) containing \(f\), after restriction near \(f\), is
comparable to a pullback of \(\hat f\), with the two induced maps to the target
differing only by a bounded modification.  This is the tropical input needed in~\cite{evc} to build local
models for the moduli stack, construct a Kuranishi structure, and hence define
Gromov--Witten invariants~\cite{vfc} and prove the associated tropical gluing
formulae~\cite{gfgw}.

\smallskip
For ordinary smooth symplectic manifolds, universal tropical structures are so easy to come by that they are hidden in the usual constructions of Gromov--Witten invariants~\cite{Ruanvirtual,FO,Tian-Li,Mcduff,pardon,GWpolyfolds,Siebert,Liu-Tian}. Indeed, if \(f\) is a nodal curve in a smooth 
manifold \(M\), then the universal tropical extension is just the orthant of
smoothing parameters for the nodes, and a family \(\hat f\) with universal
tropical structure is simply a family of curves in \(M\) containing all possible
smoothings of the nodes. In exploded and logarithmic settings~\cite{elc}, however, universal tropical structures become important and nontrivial to construct.  In logarithmic Gromov--Witten theory, universal tropical structures parallel  the `basic' condition of  Gross and Siebert~\cite{GSlogGW}, as well as the similar  `minimal log structure' condition used by Abramovich--Chen~\cite{acgw} and Kim~\cite{kim}. This condition is essential in the proof~\cite{GSlogGW} that the moduli stack
of stable logarithmic maps is algebraic: it selects a universal log structure on
the base, eliminating the excess choices present in arbitrary families of log
maps.

The importance of universal tropical structures is amplified by a useful feature of exploded
geometry.  As in logarithmic geometry, node and bubble formation can occur
within smooth families.  This makes it possible to define a natural moduli
stack of holomorphic curves, and to place it inside a natural
infinite-dimensional ambient moduli stack of not-necessarily-holomorphic
curves.  Statements formulated in this language are independent of the
particular analytic framework used to prove them.

Families of curves with universal tropical structure are used in~\cite{evc} to
prove a key structural result for this ambient stack: locally, the ambient
moduli stack retracts onto a simple finite-dimensional stack defined by a
``core family'' \(\hat f\) together with a group \(G\) of automorphisms.  In
particular, any family of curves sufficiently close to \(\hat f\) can be
written uniquely as a fiberwise holomorphic map to the domain of \(\hat f/G\),
composed with \(\hat f\), and finally modified by the exponentiation of a canonical
vector field.

Although stack language is sometimes treated with caution in symplectic
geometry, the construction of core families gives a concrete and tractable
local description of the relevant infinite-dimensional stack.  Its practical
consequence is that the \(\dbar\)-equation may be analyzed on a fixed family of
domains.  This differs in presentation from approaches based on polyfolds or
Kuranishi charts, but it is closely related to them: to construct a polyfold or Kuranishi
chart around a curve, one first needs an appropriate core family.  For curves
in smooth manifolds this preliminary step is straightforward.  In the general
exploded setting it is not.  The universal tropical structures constructed in
this paper provide the missing local models, allowing the analytic properties
of the \(\dbar\)-equation on the moduli stack to follow from the analysis
in~\cite{reg}.

\section{Notation and definitions}

The definitions in~\cite{iec} shall be essential for understanding this paper.
We recall here only the notation needed for the construction below.  An
exploded manifold \(\ex B\) has natural surjective maps of sets
\[
        \ex B \longrightarrow \totl{\ex B},
        \qquad
        \ex B \longrightarrow \totb{\ex B},
\]
to its smooth part \(\totl{\ex B}\) and tropical part \(\totb{\ex B}\);
despite the terminology, these are not subsets of \(\ex B\).  The smooth part
\(\totl{\ex B}\) is the topological space used in the sheaf-theoretic
description of \(\ex B\), analogous to the underlying scheme of a log scheme,
and should not be confused with the stronger topology in which ordinary
differential geometry on exploded manifolds takes place.

We write \(\ex T:=\et 1{\mathbb R}\) for the exploded complex torus. As a set, $\ex T$ is $\mathbb C^*\times\mathbb R$, written $\mathbb C^*\e{\mathbb R}$. 
From the perspective of differential geometry, \(\ex T\) is an infinite
disjoint union of cylinders \(\mathbb C^*\), indexed by $\mathbb R$, whereas \(\totl{\ex T}\) is a
single point, and $\totb{\ex T}$ is the line $\mathbb R$.  More generally, exploded coordinates take values in
\(\mathbb C^*\e{\mathbb R}\); thus an expression \(c\e a\) records both a
nonzero complex coefficient \(c\) and a tropical exponent \(a\).  The notation
\(\et mP\subset \ex T^m\) denotes the standard chart whose tropical part is \(P\subset\mathbb R^m\).

The tropical part \(\totb{\ex B}\) may be thought of as an integral-affine
polyhedral object, built from polytopes glued along faces.  In this paper,
however, we primarily use the  tropical structure
\((\ex B_T,\mathcal P)\), recalled below.  This records the local
integral-affine polytope over each point of \(\totl{\ex B}\), together with the
integral-affine maps obtained by moving between strata.

 We shall study the moduli stack \(\Msw(\hat{\ex B})\) of stable (not-necessarily-holomorphic) curves in a
family \(\hat{\ex B}\longrightarrow \ex G\) of exploded manifolds~\cite[Section 11]{iec}.  We shall use
the notation
\[
\begin{tikzcd}
\ex C\rar{f}\dar&\hat {\ex B}\dar \\
* \rar& \ex G
\end{tikzcd}
\]
for a \(\C\infty1\) curve \(f\) in \(\hat{\ex B}\) \cite[Definition 8.3]{iec}. Such a curve \(f\) has a tropical part \(\totb f\), a tropical curve in \(\totb{\hat{\ex B}}\), and an induced map of smooth parts \(\totl f:\totl{\ex C}\to\totl{\hat{\ex B}}\), which is a possibly nodal curve with distinguished marked points. Nodes of \(\totl{\ex C}\) correspond to internal edges of \(\totb f\), and marked points of \(\totl{\ex C}\) correspond to external edges of \(\totb f\).   We shall call \(f\) stable
if both \(f\) and the induced map of smooth parts \(\totl f\) have finite
automorphism group.  For a family \(\hat f\) of curves, we shall use the
notation
\[
\begin{tikzcd}
\ex C(\hat f)\rar{\hat f}\dar&\hat{\ex B}\dar \\
\ex F(\hat f) \rar& \ex G .
\end{tikzcd}
\]

\subsection{Tropical structure}

\

We now recall more precisely the tropical structure
\((\ex B_T,\mathcal P)\) used in this paper; see~\cite[Section 4]{iec}.
For each point \(p\in\totl{\ex B}\), the tropical structure assigns an
integral-affine polytope \(\mathcal P(p)\).  This polytope is the tropical part
of any coordinate chart whose interior stratum contains \(p\).

The polytopes \(\mathcal P(p)\) are related by canonical integral-affine face
maps.  Let
\[
        \gamma:[0,1]\longrightarrow \totl{\ex B}
\]
be a continuous path such that, whenever \(t_1\geq t_0\), the point
\(\gamma(t_1)\) lies in the closure of the stratum containing \(\gamma(t_0)\).
Then parallel transport along \(\gamma\) gives an integral-affine map
\[
        \mathcal P(\gamma):
        \mathcal P(\gamma(0))\longrightarrow \mathcal P(\gamma(1)),
\]
identifying \(\mathcal P(\gamma(0))\) with a face of
\(\mathcal P(\gamma(1))\).

For readers familiar with logarithmic geometry, this tropical structure is the
polytope-dual analogue of the ghost sheaf.  If \((X,\mathcal M_X)\) is a log
scheme, its ghost sheaf is $    \bar{\mathcal M}_X:=\mathcal M_X/\mathcal O_X^* $.
In local exploded/log models, the stalk of this ghost sheaf is the monoid of
integral-affine functions on the corresponding local tropical polytope~\cite[Section 6]{elc}.  Thus
the maps \(\mathcal P(\gamma)\) above go in the opposite direction to the
restriction maps of the ghost sheaf: restricting an integral-affine function to
a face is dual to including that face as a polytope.  This is why the tropical
structure is encoded here using polytopes and face inclusions, rather than
using monoids and restriction maps.

This data is formalized as follows.  The category \(\ex B_T\) has as objects
the points of \(\totl{\ex B}\), and as morphisms the homotopy classes of paths
\(\gamma\) satisfying the condition above.  The assignment
\[
        p\longmapsto \mathcal P(p), \qquad
        \gamma\longmapsto \mathcal P(\gamma)
\]
defines a functor from \(\ex B_T\) to the category of integral-affine
polytopes.  We denote this functor again by \(\mathcal P\), and call the pair
\((\ex B_T,\mathcal P)\) the tropical structure of \(\ex B\).

This construction is functorial.  Given a map of exploded manifolds
\[
        h:\ex A\longrightarrow \ex B ,
\]
there is an induced map of tropical structures
\[
        (h_T,\mathcal P h):
        (\ex A_T,\mathcal P_{\ex A})
        \longrightarrow
        (\ex B_T,\mathcal P_{\ex B}).
\]
Here \(h_T:\ex A_T\longrightarrow \ex B_T\) is the functor induced by the map
of smooth parts
\[
        \totl h:\totl{\ex A}\longrightarrow \totl{\ex B},
\]
sending a path \(\gamma\) to \(\totl h\circ\gamma\).  The second component is a
natural transformation
\[
        \mathcal P h:
        \mathcal P_{\ex A}
        \longrightarrow
        \mathcal P_{\ex B}\circ h_T .
\]
At a point \(x\in\totl{\ex A}\), this is the integral-affine map
\[
        \mathcal P h:
        \mathcal P_{\ex A}(x)
        \longrightarrow
        \mathcal P_{\ex B}(\totl h(x))
\]
given by the tropical part of \(h\) in coordinate charts whose interior strata
contain \(x\) and \(\totl h(x)\).

When no confusion is possible, we suppress the subscripts and write
\(\mathcal P\) for the polytope functor associated to any exploded manifold.
Throughout this paper, an arrow
\[
        (\ex A_T,\mathcal P)\longrightarrow(\ex B_T,\mathcal P)
\]
will always mean a functor on the categories of strata together with the
corresponding natural transformation of polytope functors, as above.

\subsection{Extensions of the tropical structure of a curve}

\

The purpose of this subsection is to define the tropical structures which can
arise by placing a fixed curve \(f\) inside a family of curves.  Such a family
may extend the tropical
data over \(f\): the one-point tropical base of \(f\) is replaced by a polytope
of possible tropical deformation parameters.  The resulting object is what we
call an extended tropical structure on \(f_T\).

The tropical structure of a family of curves \(\hat f\) in \(\hat{\ex B}\) is
the diagram
\[
\begin{tikzcd}
(\ex C(\hat f)_T,\mathcal P)\rar{(\hat f_T,\mathcal P\hat f)}\dar
    &(\hat{\ex B}_T,\mathcal P)\dar \\
(\ex F(\hat f)_T,\mathcal P)\rar
    &(\ex G_T,\mathcal P).
\end{tikzcd}
\]
In the case of a single curve \(f\), the base \(\ex F(f)\) is a single point
\(*\), and the associated polytope is also a single point.  We shall therefore
write the tropical structure of \(f\) as
\[
\begin{tikzcd}
(\ex C_T,\mathcal P)\rar{(f_T,\mathcal P f)}\dar
    &(\hat{\ex B}_T,\mathcal P)\dar \\
(*,*)\rar
    &(\ex G_T,\mathcal P).
\end{tikzcd}
\]
We denote this tropical structure of $f$ by \(\mathcal P_f\). 

Now suppose that a family \(\hat f\) contains \(f\).  Then the inclusion of
\(f\) as a fibre of \(\hat f\) gives a commutative diagram
\[
\begin{tikzcd}
(\ex C_T,\mathcal P)\arrow{r}\dar\ar[bend left]{rr}{(f_T,\mathcal P f)}
    &(\ex C(\hat f)_T,\mathcal P)\rar\dar
    &(\hat{\ex B}_T,\mathcal P)\dar \\
(*,*)\rar
    &(\ex F(\hat f)_T,\mathcal P)\rar
    &(\ex G_T,\mathcal P).
\end{tikzcd}
\]
Restricting the tropical structure of the total space \(\ex C(\hat f)\) to the
image of \(\ex C\) gives a new functor, which we denote by \(\mathcal P'\), on
\(\ex C_T\).  Similarly, restricting the tropical structure of the base
\(\ex F(\hat f)\) to the image of the point \(*\) gives a polytope
\(\mathcal P'(*):=P\), together with a map
\[
(*,P)\longrightarrow(\ex G_T,\mathcal P).
\]
Here, we use the following shorthand: if \(P\) is an integral-affine polytope,
then \((*,P)\) denotes the tropical structure on the one-object category \(*\)
whose polytope functor takes the unique object to \(P\).

Thus the family \(\hat f\) determines a diagram
\[
\begin{tikzcd}
(\ex C_T,\mathcal P')\rar\dar
    &(\hat{\ex B}_T,\mathcal P)\dar \\
(*,P)\rar
    &(\ex G_T,\mathcal P).
\end{tikzcd}
\]
This is the basic example of an extended tropical structure on \(f_T\).

The inclusion of \(f\) as a fibre of \(\hat f\) gives more than an extended tropical structure: it gives a distinguished morphism from the original tropical structure \(\mathcal P_f\) to this extended tropical structure. Equivalently, it chooses the point of the deformation polytope \(P\) corresponding to the original curve \(f\), together with an identification of the pullback over that point with \(\mathcal P_f\). This is expressed by the left-hand pullback square in the diagram
\[
\begin{tikzcd}
(\ex C_T,\mathcal P)\rar\dar
    &(\ex C_T,\mathcal P')\rar\dar
    &(\hat{\ex B}_T,\mathcal P)\dar \\
(*,*)\rar
    &(*,P)\rar
    &(\ex G_T,\mathcal P).
\end{tikzcd}
\]
This distinguished morphism is the data which will later be called an extension of the tropical structure of \(f\). This notion is the polytope-dual analogue of Gross--Siebert's ghost sheaf of a given type~\cite{GSlogGW}. In Theorem~\ref{mts}, we construct a universal such extension.

 \begin{defn}\label{ets}
An \emph{extended tropical structure} on \(f_T\) consists of a base polytope
\(P\), a functor
\[
        \mathcal P':\ex C_T\longrightarrow
        \{\text{integral-affine polytopes}\},
\]
a map \((*,P)\longrightarrow(\ex G_T,\mathcal P)\), and a commutative diagram
\begin{equation}\label{etscd}
\begin{tikzcd}[column sep=large]
(\ex C_T,\mathcal P')\dar{(\pi_T,\mathcal P'\pi)}
        \rar{(f_T,\mathcal P'f)}
    &(\hat{\ex B}_T,\mathcal P)\dar \\
(*,P)\rar
    &(\ex G_T,\mathcal P),
\end{tikzcd}
\end{equation}
satisfying the following conditions.
\begin{enumerate}
\item
For every morphism \(\gamma\) of \(\ex C_T\),
\[
        \mathcal P'(\gamma):
        \mathcal P'(\gamma(0))\longrightarrow \mathcal P'(\gamma(1))
\]
is an integral-affine isomorphism onto a face of
\(\mathcal P'(\gamma(1))\).

\item\label{ets2}
For every point \(x\in\totl{\ex C}\), the map
\[
        \mathcal P'\pi:\mathcal P'(x)\longrightarrow P
\]
has the following local form.
\begin{enumerate}
\item If \(x\) is in a smooth component of \(\ex C\), then
\(\mathcal P'\pi\) is an isomorphism.

\item If \(x\) is a node of \(\totl{\ex C}\), then \(\mathcal P'\pi\) is a
pullback of the map
\[
        [0,\infty)^2\longrightarrow [0,\infty),
        \qquad
        (a,b)\longmapsto a+b;
\]
that is, there is a pullback diagram
\[
\begin{tikzcd}
\mathcal P'(x)\rar\dar
    & {[}0,\infty)^2\dar{a+b} \\
P\rar
    & {[}0,\infty).
\end{tikzcd}
\]

\item If \(x\) is a marked point of \(\totl{\ex C}\), equivalently if the
inverse image of \(x\) in \(\ex C\) is isomorphic to
\(\et 1{(0,\infty)}\), then \(\mathcal P'\pi\) is a pullback of the map
\[
        [0,\infty)\longrightarrow\{0\};
\]
that is, there is a pullback diagram
\[
\begin{tikzcd}
\mathcal P'(x)\rar\dar
    & {[}0,\infty)\dar \\
P\rar
    & \{0\}.
\end{tikzcd}
\]
\end{enumerate}
\end{enumerate}
\end{defn}

Condition~\ref{ets2} is the tropical form of the local models for a family of
curves: over a smooth point the local tropical polytope is just the base
polytope, at a node there are two branch parameters whose sum is the smoothing
parameter, and at a marked point there is one additional unbounded direction.

Provided that \(\ex C\neq\ex T\), any family \(\hat f\) of curves containing
\(f\) gives an extended tropical structure on \(f_T\) by the restriction
procedure described above.  We shall ignore the case \(\ex C=\ex T\), since the
moduli stack of curves with this domain is elementary to analyze.

\
The correct notion of changing the base polytope of an extended tropical
structure is pullback.  If \(\mathcal P'\) is an extended tropical structure
with base polytope \(P\), and if
\[
        \phi:Q\longrightarrow P
\]
is an integral-affine map of polytopes, define
\[
        \phi^*\mathcal P'(x):=Q\times_P\mathcal P'(x)
\]
for every \(x\in\totl{\ex C}\).  The face maps and the maps to
\((\hat{\ex B}_T,\mathcal P)\) are induced by the corresponding maps for
\(\mathcal P'\).

 \begin{lemma}\label{pets}
The construction above defines an extended tropical structure
\(\phi^*\mathcal P'\) on \(f_T\) with base polytope \(Q\).  It is the unique
extended tropical structure fitting into a pullback diagram
\[
\begin{tikzcd}[column sep=large]
(\ex C_T,\phi^*\mathcal P')\dar
        \ar[bend left]{rr}{(f_T,\phi^*\mathcal P'f)}
        \rar
    &(\ex C_T,\mathcal P')\rar{(f_T,\mathcal P'f)}\dar
    &(\hat{\ex B}_T,\mathcal P)\dar \\
(*,Q)\rar
    &(*,P)\rar
    &(\ex G_T,\mathcal P),
\end{tikzcd}
\]
in the sense that, for every \(x\in\totl{\ex C}\), the square
\[
\begin{tikzcd}
\phi^*\mathcal P'(x)\rar\dar
    &\mathcal P'(x)\dar \\
Q\rar{\phi}
    &P
\end{tikzcd}
\]
is a pullback square.
\end{lemma}
 
 \begin{proof}
The functoriality of \(\phi^*\mathcal P'\) follows from the universal property
of the pullback squares defining \(\phi^*\mathcal P'(x)\).  Since face
inclusions, the node model
\([0,\infty)^2\to[0,\infty)\), and the marked-point model
\([0,\infty)\to\{0\}\) are all preserved by base change, the two conditions in
Definition~\ref{ets} hold for \(\phi^*\mathcal P'\).  The same pullback
universal property gives uniqueness.
\end{proof}

 Thus, if \(\mathcal Q\) and \(\mathcal P'\) are extended tropical structures with base polytopes \(Q\) and \(P\), respectively, a map of extended tropical structures from \(\mathcal Q\) to \(\mathcal P'\) will mean precisely such a pullback diagram:
\[
\begin{tikzcd}
(\ex C_T,\mathcal Q)\ar[bend left]{rr}\rar\dar &(\ex C_T,\mathcal P')\rar\dar &(\hat{\ex B}_T,\mathcal P)\dar \\
(*,Q)\rar &(*,P)\rar &(\ex G_T,\mathcal P).
\end{tikzcd}
\]

In particular, the tropical structure \(\mathcal P_f\) of $f$ is an extended tropical structure with base polytope a point.
 
 \begin{defn}
An \emph{extension of the tropical structure of \(f\)} is an extended tropical structure \(\mathcal P'\) on \(f_T\), with base polytope \(P\), together with a chosen morphism
\[
        \iota:\mathcal P_f\longrightarrow \mathcal P'.
\]
Equivalently, it is a diagram
\[
\begin{tikzcd}
(\ex C_T,\mathcal P)\rar\dar &(\ex C_T,\mathcal P')\rar\dar &(\hat{\ex B}_T,\mathcal P)\dar \\
(*,*)\rar &(*,P)\rar &(\ex G_T,\mathcal P)
\end{tikzcd}
\]
whose left-hand square is a pullback. The chosen map \(*\to P\) is the point of the deformation polytope corresponding to the original curve \(f\).
\end{defn}

%

\begin{defn} \label{universal extension}
A \emph{universal extension} of the tropical structure of \(f\) is an extension
\[
        \iota_u:\mathcal P_f\longrightarrow \mathcal P_u
\]
with the following property: for any extension
\[
        \mathcal P_f\longrightarrow \mathcal P',
\]
there is a unique map of extended tropical structures
\[
        \mathcal P'\longrightarrow \mathcal P_u
\]
such that the resulting triangle commutes. 
\[\begin{tikzcd} \mathcal P_f \rar \dar{\iota_u} & \mathcal P'\ar[dotted]{dl}{\exists!}
\\ \mathcal P_u\end{tikzcd}\]
 Equivalently, every extension of
the tropical structure of \(f\) is obtained uniquely as a pullback of the
universal extension, compatibly with the chosen morphism from the tropical
structure of \(f\).
\end{defn}

 \section{Construction of a universal extension of tropical structure}
 
 \begin{thm}\label{mts}
Let \(f\) be a curve in \(\hat{\ex B}\) with domain \(\ex C\neq\ex T\).  Then there exists an extension
\[
        \iota_u:\mathcal P_f\longrightarrow \mathcal P_u
\]
which is universal in the sense of Definition~\ref{universal extension}.
\end{thm}

 \begin{proof}

Let \(P\) denote the base polytope of an arbitrary extension \(\iota':\mathcal P_f\to\mathcal P'\), and let \(P_u\) denote the base polytope to be constructed for the universal extension \(\mathcal P_u\).

The construction has four parts.
\begin{enumerate}
\item Construct an ambient parameter polytope \(Q\) recording unconstrained vertex-position and edge-length data. Over \(Q\), construct the domain-side data \(\mathcal Q\) of a putative extended tropical structure: a functor from \(\ex C_T\) to integral-affine polytopes, together with maps to \(Q\) having the required local form over smooth points, nodes, and marked points.
\item Show that any extension \(\iota':\mathcal P_f\to\mathcal P'\) determines a unique map \(P\to Q\), and that \(\mathcal P'\) is the pullback of this domain-side data.
\item Identify the integral-affine consistency equations required for the candidate maps \(\mathcal Q(x)\to\mathcal P(f(x))\), initially defined on smooth strata, to be independent of choices and to extend over marked points and nodes as a natural transformation.
\item Define \(P_u\subset Q\) by these equations. Pulling \(\mathcal Q\) back to \(P_u\) gives an extended tropical structure \(\mathcal P_u\), and the chosen point corresponding to \(f\) gives an extension \(\iota_u:\mathcal P_f\to\mathcal P_u\). The uniqueness in the previous steps then gives the universal property.
\end{enumerate}

\smallskip
\noindent\emph{Step 1: the unconstrained domain-side data.}

Choose a point \(p_i\) on each smooth component of \(\ex C\). For each node \(e\) of \(\totl{\ex C}\), corresponding to an internal edge of \(\totb{\ex C}\), set \(I_e:=[0,\infty)\). Define
\begin{equation}\label{Qdef}Q:=\prod_{i}\mathcal P(f(p_{i}))\times \prod_{e}I_{e}\ . \end{equation}
The factors \(\mathcal P(f(p_i))\) record possible positions of the vertices of the tropical curve \(\totb f\), while the factors \(I_e\) record possible lengths of its internal edges. At this stage \(Q\) is only an unconstrained parameter polytope. The functor \(\mathcal Q\) constructed below has the correct domain-side local form over \(Q\), but it is not yet an extended tropical structure, because the maps to the tropical structure of \(\hat{\ex B}\) have not yet been made into a natural transformation. The universal deformation polytope \(P_u\subset Q\) will be the subpolytope cut out by  the consistency equations needed for this natural transformation to exist.

 \smallskip
 
 As a first step to define a functor $\mathcal Q$ from $\ex C_{T}$ to the category of integral-affine polytopes, for each point $x$ on a smooth component of $\ex C$, define \[\mathcal Q(x):=Q\] For any path $\gamma$ within a smooth component, define the map $\mathcal Q(\gamma)$ to be the identity.
 
  For any $x$ in a smooth component of $\ex C$, define 
 \[\mathcal Qf:\mathcal Q(x)\longrightarrow \mathcal P(f(x))\]
by choosing a path $\gamma$ within a stratum of $\totl{\ex C}$ joining $p_{i}$ to $x$, and setting $\mathcal Qf$ to be projection to $\mathcal P(f(p_{i}))$ followed by $\mathcal P(f_{T}(\gamma))$. Note that $\mathcal Qf$ may depend on this choice of path, in which case it won't be a natural transformation. Later, we shall restrict to subpolytopes so that we get a natural transformation.

\smallskip

For any node $e$ of $\totl{\ex C}$ (corresponding to an internal edge of $\totb{\ex C}$), define $\mathcal Q(e)$ to be the product of $[0,\infty)^{2}$ with $I_{e'}$ for every other node $e'$, and $\mathcal P(f(p_{i}))$ for all $i$.
\[\mathcal Q(e):=[0,\infty)^{2}\times \prod_{i}\mathcal P(f(p_{i}))\times \prod_{e'\neq e}I_{e'}\]
 For a path $\gamma$ joining $x$ to one side of the node $e$, set $\mathcal Q(\gamma)$ to be the map induced by inclusion of  $I_{e}:=[0,\infty)$ as  $\{0\}\times[0,\infty)\subset [0,\infty)^{2}$. (Let $\mathcal Q(\gamma)$ be the identity on all other factors.) For a path $\gamma'$ joining $x'$ to the other side of the node $e$, let $\mathcal Q(\gamma')$ be the inclusion of  $I_{e}$ as  $[0,\infty)\times\{0\}\subset [0,\infty)^{2}$ and the identity on all other factors of $\mathcal Q(x')=Q$.

\smallskip

A marked point $y$ of $\totl{\ex C}$ is a point in $\totl{\ex C}$ with inverse image in $\ex C$ isomorphic to $\et 1{(0,\infty)}$, so $y$ corresponds to an external edge of $\totb{\ex C}$. For any marked point $y$ of $\totl{\ex C}$, define $\mathcal Q(y)=Q\times [0,\infty)$. For a path $\gamma$ joining a point $x$ in a smooth component  of $\ex C$ to $y$, define $\mathcal Q(\gamma):\mathcal Q(x)\longrightarrow\mathcal Q(y)$ to be the inclusion of $Q$ as the face $Q\times \{0\}\subset Q\times[0,\infty)$.

\smallskip

This completes the definition of a functor $\mathcal Q$ from $\ex C_{T}$ to the category of integral-affine polytopes. Define $\mathcal Q(*)$ to be $Q$. There is a map 
\[(\ex C_{T},\mathcal Q)\longrightarrow (*, Q)\]
given by
\begin{itemize}\item the identity map \[\mathcal Q(x)= Q\longrightarrow Q\] when $x$ is in a smooth component of $\ex C$, \item the projection
\[\mathcal Q(y):=Q\times[0,\infty)\longrightarrow Q\]
when  $y$ is a marked point of $\totl{\ex C}$,
\item and the map
\[[0,\infty)^{2}\longrightarrow I_{e}\]
\[(a,b)\mapsto a+b\]
on the components of $\mathcal Q(e)$ and $\mathcal Q(*)$ corresponding to the node $e$, and the identity map on all other factors of $\mathcal Q(e)$ and $\mathcal Q(*)$.
\end{itemize}

\smallskip
\noindent\emph{Step 2: the universal property before imposing consistency equations.}

We have not yet defined a natural transformation \(\mathcal Qf\), only the candidate maps \(\mathcal Q(x)\to\mathcal P(f(x))\) for \(x\) in smooth components of \(\ex C\). Nevertheless, \(\mathcal Q\) already represents the unconstrained domain-side data of any extension.

\begin{claim}\label{unique inclusion}
Given any extension \(\mathcal P_f\to\mathcal P'\) with base polytope \(P\), there exists a unique natural transformation
\[
        \eta_0:\mathcal P'\longrightarrow\mathcal Q
\]
so that
\begin{enumerate}\item\label{ui1} the following diagram commutes
\[\begin{tikzcd}(\ex C_{T},\mathcal P')\dar \rar{(\id,\eta_{0})}&(\ex C_{T},\mathcal Q)\dar
\\(*, P)\rar &(*, Q)
\end{tikzcd}\]
\item \label{ui2}for all $x$ in a smooth component of  $\ex C$,
\[\mathcal P'f:\mathcal P'(x)\longrightarrow \mathcal P(f(x))\]
factorizes as 
\[\mathcal P'f=\mathcal Qf\circ\eta_{0}\]
\item \label{ui3}and for all points $x$ in $\totl{\ex C}$, the following is a pullback diagram:
\[\begin{tikzcd}\mathcal P'(x)\dar\rar{\eta_{0}}&\mathcal Q(x)\dar
\\ P\rar &Q\end{tikzcd}\]
\end{enumerate}\end{claim}

To prove Claim \ref{unique inclusion}, note that condition \ref{ets2} of Definition \ref{ets} tells us that at points $x$ in the smooth part of $\ex C$, the map $\mathcal P'(x)\longrightarrow P$ must be an isomorphism. The commutativity of the diagram in item \ref{ui1} above and the fact that $\mathcal Q(x)\longrightarrow Q$ is an isomorphism tells us that $\mathcal P'(x)\longrightarrow \mathcal Q(x)$ must be determined by the map $P\longrightarrow Q$. On the other hand, item \ref{ui2} above tells us that the map $\mathcal P'(p_{i})=P\longrightarrow Q=\mathcal Q(p_{i})$ followed by projection to $\mathcal P(f(p_{i}))$ must be equal to $\mathcal P'f$.  

Condition \ref{ets2} of Definition \ref{ets} and the definition of a pullback in Lemma \ref{pets} imply  that there is a pullback diagram:
\[\begin{tikzcd}\mathcal P'(e)\rar \dar&\mathcal Q(e)\rar\dar &{} [0,\infty)^{2} \dar 
\\ P \rar& Q\rar&{} [0,\infty)\end{tikzcd}\]

The map $P\longrightarrow[0,\infty)$ coming from the node $e $ must therefore be equal to the map $P\longrightarrow Q$ followed by projection to  $I_{e}=[0,\infty)$. The map $P\longrightarrow Q$ is therefore uniquely determined by the above conditions. This implies that $\eta_{0}:\mathcal P'(x)\longrightarrow \mathcal Q(x) $ is uniquely determined for any point $x$ in a smooth component of $\ex C$.

The fact that $\mathcal P$ is a pullback of $\mathcal P'$ along with condition \ref{ets2} of Definition \ref{ets} implies that given paths $\gamma$ and $\gamma'$ ending on opposite sides of the node $e$, the maps $\mathcal P'(\gamma)$ and $\mathcal P'(\gamma')$ send $P$ to the two different faces of $\mathcal P'(e)$ which are the pullback of the two different  boundary faces of $[0,\infty)^{2}$. Therefore, there is a unique map  $\eta_{0}:\mathcal P'(e)\longrightarrow \mathcal Q(e)$ so that item \ref{ui3} above holds and the following diagram commutes, as is required if $\eta_{0}$ is to be a natural transformation.
\[\begin{tikzcd}
 P\rar{\eta_{0}}\dar[swap]{\mathcal P'(\gamma')} &Q\dar{\mathcal Q(\gamma')}
\\\mathcal P'(e)\rar{\eta_{0}}&\mathcal Q(e)
\\ P\rar{\eta_{0}}\ar{u}{\mathcal P'(\gamma)} &Q\uar[swap]{\mathcal Q(\gamma)}\end{tikzcd}\]

Similarly, for a marked point $y$ of $\totl{\ex C}$, \[\eta_{0}:\mathcal P'(y)\longrightarrow\mathcal Q(y)\] is uniquely determined by condition \ref{ets2} of Definition \ref{ets} and item \ref{ui3} above. Given any path $\gamma$ joining $x$ to $y$, the diagram required to show that $\eta_{0}$ is a natural transformation commutes
\[\begin{tikzcd}\mathcal P'(y)\rar{\eta_{0}}&\mathcal Q(y)
\\\mathcal P'(x)\rar{\eta_{0}}\uar{\mathcal P'(\gamma)}&\mathcal Q(x)\uar{\mathcal Q(\gamma)}\end{tikzcd}\]
because in each case the inclusion corresponding to $\gamma$ is the pullback of the inclusion ${0}\longrightarrow[0,\infty)$.

It follows that $\eta_{0}:\mathcal P'\longrightarrow\mathcal Q$ as defined is a natural transformation, and is the unique natural transformation satisfying the above conditions. This completes the proof of Claim \ref{unique inclusion}.

\

\smallskip
\noindent\emph{Step 3: the consistency equations.}

So far, the maps \(\mathcal Qf:\mathcal Q(x)\to\mathcal P(f(x))\) have only been defined for \(x\) in smooth components of \(\ex C\), and they may depend on choices of paths. We now impose the integral-affine conditions needed for these maps to be independent of choices and to extend over marked points and nodes as a natural transformation.

  For $\mathcal Qf$ to be a natural transformation, we need for any path $\gamma$ joining $x_{1}$ to $x_{2}$ within a given  smooth component, the following diagram to commute:
 \[\begin{tikzcd}Q=\mathcal Q(x_{1})\rar{\mathcal Qf}\dar{\id=\mathcal Q(\gamma)}&\mathcal P(f(x_{1}))\dar{\mathcal P(f_{T}(\gamma))}
 \\  Q=\mathcal Q(x_{2})\rar{\mathcal Qf}&\mathcal P(f(x_{2}))\end{tikzcd}\]
If $x_{1}=x_{2}$, then this is not possible unless $\mathcal Qf$ has image contained in the subset of $\mathcal P(f(x_{1}))$ where $\mathcal P(f_{T}(\gamma))$ is the identity.

 Define \[\mathcal P_{1}(f(x))\subset \mathcal P(f(x))\] to be the subpolytope of $\mathcal P(f(x))$ which is fixed by $\mathcal P(f_{T}(\gamma))$ for all such $\gamma$ from $x$ to itself. Define $Q_{1}$ to be the subpolytope of $Q$ which is the product of $\mathcal P_{1}(f(p_{i}))$ for all $i$ with $I_{e}$ for each node $e$. 
 
 \[Q_{1}:=\prod_{i}\mathcal P_{1}(f(p_{i}))\times \prod_{e}I_{e}\]
 
 As $\mathcal P'f=\mathcal Qf\circ \eta_0$ restricted to $x$, and $\mathcal P'(\gamma)$ is the identity for any loop $\gamma$, the map $P\longrightarrow Q$ from Claim \ref{unique inclusion} must have image inside $Q_{1}$.
%
%
 
 \
 
 Now define \[\mathcal Qf:\mathcal Q(y)\longrightarrow \mathcal P(f(y))\] for a marked point $y$ of $\totl{\ex C}$. Choose a path $\gamma$ to $y$ from a point $x$ in the smooth component of $\totl{\ex C}$ with closure containing $y$.
For $\mathcal Qf$ to be a natural transformation, we want the following diagram to commute.
\begin{equation}\label{ym}\begin{tikzcd}\mathcal Q(y)\rar{\mathcal Qf}&\mathcal P(f(y))
\\ \mathcal Q(x)\uar{\mathcal Q(\gamma)}\rar{\mathcal Qf}&\uar{\mathcal P(f_{T}(\gamma))}\mathcal P(f(x))\end{tikzcd}\end{equation}
 As $\mathcal Q(y)=Q\times[0,\infty)$ and $\mathcal Q(\gamma)$ has image $Q\times\{0\}$, this specifies $\mathcal Qf$ on the face $Q\times\{0\}\subset Q\times[0,\infty)$. Claim \ref{unique inclusion}  gives a unique inclusion $\eta_{0}:\mathcal P(y)\longrightarrow \mathcal Q(y)$ as some fiber of the projection $Q\times[0,\infty)\longrightarrow Q$. As the image of $\mathcal P(y)$ is complementary to $Q\times\{0\}$, we may define a unique map 
 \[\mathcal Qf:\mathcal Q(y)\longrightarrow \mathcal P(f(y))\]
 so that diagram (\ref{ym}) commutes and so that $\mathcal Pf= \mathcal Qf\circ \eta_{0}$. Note that $\mathcal Qf$ as defined depends on the choice of $\gamma$, but restricted to the inverse image of $Q_{1}\subset Q$, it does not depend on the choice of $\gamma$.
 
 The map $\eta_{0}:\mathcal P'(y)\longrightarrow \mathcal Q(y)$ pulls back $\mathcal Qf$ to $\mathcal P'f$ because the following diagram commutes.
\[\begin{tikzcd}\mathcal P(y)\rar{\eta_{0}} \ar[bend left]{rr}{\mathcal Pf}\ar{dr}&\mathcal Q(y)\rar{\mathcal Qf}&\mathcal P(f(y))  
\\ &\mathcal P'(y)\uar{\eta_{0}} \ar{ur}{\mathcal P'f}\end{tikzcd}\]
In particular, the bottom left hand triangle commutes because of the uniqueness property of $\eta_{0}$, the top loop commutes by the definition of $\mathcal Qf$, and the outer loop commutes because $\mathcal P$ is a pullback of $\mathcal P'$. It follows that the bottom right hand triangle commutes on the image of $\mathcal P(y)$. On the other hand, the bottom right hand triangle commutes on the complementary  face of $\mathcal P'(y)$ because of condition \ref{ui2} from Claim \ref{unique inclusion} and the commutativity of diagram (\ref{ym}).
 
 \smallskip
 
 Let us now try to define $\mathcal Qf:\mathcal Q(e)\longrightarrow \mathcal P(f(e))$ similarly for a node $e$ of $\totl{\ex C}$. 
  Given a path $\gamma$  joining  $p_{i}$  to the node $e$, if $\mathcal Qf$ is to be a natural transformation, there must be a commutative diagram
  \[\begin{tikzcd}\mathcal Q(e)\rar[dotted]{\mathcal Qf?}&\mathcal P(f(e))
  \\\mathcal Q(p_{i})\uar{\mathcal Q(\gamma)}\rar{\mathcal Qf}&\mathcal P(f(p_{i}))\uar{\mathcal P(f_{T}(\gamma))}
  \end{tikzcd}\]
  
Claim \ref{unique inclusion} gives us an inclusion of $\mathcal P(e)$ as a fiber of $\mathcal Q(e)\longrightarrow Q$. As the pullback of $\mathcal Qf$ under this inclusion is required to be $\mathcal Pf$, we now have our prospective map $\mathcal Q f$ specified on the face which is the image of $\mathcal Q(\gamma)$ and the complementary subpolytope $\mathcal P(e)\subset\mathcal Q(e)$. If we embed $\mathcal P(f(e))$ in $\mathbb R^{n}$, this specifies a unique integral-affine map
\[A_{\gamma}:\mathcal Q(e)\longrightarrow \mathbb R^{n}\supset\mathcal P(f(e))\] 
so that the following diagram commutes
\begin{equation}\label{Adef}\begin{tikzcd} 
 \mathcal P(e)\rar{\eta_{0}}\ar[bend left]{rr}{\mathcal Pf}&\mathcal Q(e) \rar{A_{\gamma}}&\mathbb R^{n}\supset\mathcal P(f(e)) 
\\ &\mathcal Q(x)\uar{\mathcal Q(\gamma)}\rar{\mathcal Qf}&\mathcal P(f(x))\uar[swap]{\mathcal P(f_{T}(\gamma))}\end{tikzcd}\end{equation}

Consider the  following diagram.
 \[\begin{tikzcd}\mathcal P(e)\rar{\eta_{0}} \ar[bend left]{rr}{\mathcal Pf}\ar{dr}&\mathcal Q(e)\rar{A_{\gamma}}&\mathbb R^{n}\supset \mathcal P(f(e))  
\\ &\mathcal P'(e)\uar{\eta_{0}} \ar{ur}{\mathcal P'f}\end{tikzcd}\]
The top, lower left and outer loops of the above diagram commute, so the bottom righthand triangle commutes when $\mathcal P'(e)$ is restricted to the image of $\mathcal P(e)$. On the other hand, consider the diagram 
\[\begin{tikzcd} 
 \mathcal P'(e)\rar{\eta_{0}} \ar[bend left]{rr}{\mathcal P'f}&\mathcal Q(e) \rar{A_{\gamma}}&\mathbb R^{n}\supset\mathcal P(f(e)) 
\\ \mathcal P'(x)\uar{\mathcal P'(\gamma)}\rar\ar[bend right]{rr}{\mathcal P'f}&\mathcal Q(x)\uar{\mathcal Q(\gamma)}\rar{\mathcal Qf}&\mathcal P(f(x))\uar[swap]{\mathcal P(f_{T}(\gamma))}\end{tikzcd}\]
The outer loop, two inner squares and bottom triangle commute, therefore the top triangle commutes when $\mathcal P'(e)$ is restricted to the image of $\mathcal P'(x)$. Therefore both of the above diagrams commute, because the questionable triangle commutes when $\mathcal P'(e)$ is restricted to the complementary  images of $\mathcal P(e)$ and $\mathcal P'(x)$.  In particular, this implies that the pullback of $A_{\gamma}$ to $\mathcal P'(e)$ must   equal  $\mathcal P'f$.

\smallskip

Given another path \(\gamma'\) joining a point \(p_j\) on the other side of the node \(e\) to \(e\), we get an analogous map 
\[A_{\gamma'}:\mathcal Q(e)\longrightarrow \mathbb R^{n}\supset\mathcal P(f(e))\]
As the diagram 
\[\begin{tikzcd}\mathcal Q(e)\rar{A_{\gamma}}&\mathbb R^{n}  
\\ \mathcal P'(e)\uar{\eta_{0}} \ar{ur}{\mathcal P'f}\rar{\eta_{0}}& \mathcal Q(e)\uar[swap]{A_{\gamma'}}\end{tikzcd}\]
  commutes,
  $\eta_{0}(\mathcal P'(e))$ must be in the subset of $\mathcal Q(e)$ on which $A_{\gamma}=A_{\gamma'}$. This subset of $\mathcal Q(e)$ is the inverse image of the subset of $ Q=\mathcal Q(p_{j})$ on which the diagram 
\[  \begin{tikzcd}\mathcal Q(p_{j})\rar{\mathcal Qf} \dar{\mathcal Q(\gamma')}&\mathcal P(f(p_{j}))
\dar{\mathcal P(f_{T}\gamma')}
\\ \mathcal Q(e)\rar{A_{\gamma}}&\mathbb R^{n}\supset \mathcal P(f(e))
 \end{tikzcd}\] commutes.
 
\

\smallskip
\noindent\emph{Step 4: construction of the universal extension.}

Define \(P_u\subset Q_1\subset Q\) to be the subpolytope on which the compatibility diagram above commutes for every node \(e\) and every pair of paths \(\gamma,\gamma'\) to \(e\) from the two sides.  Equivalently, \(P_u\) is the solution set of the integral-affine consistency equations just described: consistency equations from monodromy within smooth components cut out $Q_1\subset Q$, and then consistency equations from nodes cut out $P_u\subset Q_1$.

Define \(\mathcal P_u\) as the pullback of \(\mathcal Q\) along the inclusion \(P_u\hookrightarrow Q\). In other words define $\mathcal P_{u}(x)$ via the pullback diagram
\[\begin{tikzcd}\mathcal P_{u}(x)\rar\dar&\mathcal Q(x)\dar
\\ P_{u}\rar&Q \end{tikzcd}\]
and define $\mathcal P_{u}(\gamma)$ to be the map induced by the diagram
\[\begin{tikzcd}\mathcal P_{u}(x_{1})\rar\ar[bend left]{rr}\dar&\mathcal Q(x_{1})\ar{dr}\rar{\mathcal Q(\gamma)}&\mathcal Q(x_{2})\dar
\\ P_{u}\ar{rr}&&Q \end{tikzcd}\]

For points $x$ in $\totl{\ex C}$ which are not nodes, we may define $\mathcal P_{u}f$ via the composition
\[\begin{tikzcd}\mathcal P_{u}(x)\ar[bend left]{rr}{\mathcal P_{u}f}\rar&\mathcal Q(x)\rar{\mathcal Qf}&\mathcal P(f(x))\end{tikzcd}\]

To define $\mathcal P_{u}f$ at a node $e$ of $\totl{\ex C}$, consider the commutative diagram 
\[\begin{tikzcd}\mathcal P_{u}(p_{j})\dar[swap]{\mathcal P_{u}(\gamma')}\ar{rr}{\mathcal P_{u}f}&&\mathcal P(f(p_{j}))\dar{\mathcal P(f_{T}\gamma')}
\\ \mathcal P_{u}(e)\rar&\mathcal Q(e)\rar{A_{\gamma}} &\mathbb R^{n}\supset\mathcal P(f(e))
\\ \mathcal P_{u}(p_{i})\uar{\mathcal P_{u}(\gamma)}\ar{rr}{\mathcal P_{u}f}&&\mathcal P(f(p_{i}))\uar[swap]{\mathcal P(f_{T}\gamma)}
\end{tikzcd}\]
As $\mathcal P_{u}(e)$ is the convex hull of the image of $\mathcal P_{u}(\gamma)$ and $\mathcal P_{u}(\gamma')$, and all polytopes above are convex, the composition 
\[\begin{tikzcd}\mathcal P_{u}(e)\rar&\mathcal Q(e)\rar{A_{\gamma}} &\mathbb R^{n}\supset\mathcal P(f(e))\end{tikzcd}\]
has image inside $\mathcal P(f(e))$, so we may define $\mathcal P_{u}f$ using the commutative diagram
\[\begin{tikzcd}\mathcal P_{u}(e)\ar[bend left]{rrr}{\mathcal P_{u}f}\rar& \mathcal Q(e)\rar{A_{\gamma}}&\mathbb R^{n}\rar[hookleftarrow]& \mathcal P(f(e))\end{tikzcd}\]
As discussed above, $\mathcal P_{u}f$ thus defined obeys the commutativity requirements to be a natural transformation. 
This completes the construction of our universal extension of the tropical structure of $f$
\[ \begin{tikzcd}[column sep=large](\ex C_{T},\mathcal P_{u})\dar\rar{(f_{T},\mathcal P_{u}f)}&(\hat{\ex B}_{T},\mathcal P)\dar
\\(*, P_{u})\rar&(\ex G_{T},\mathcal P) \end{tikzcd} \]

Now let \(\iota':\mathcal P_f\to\mathcal P'\) be any extension with base polytope \(P\). By Step 2, it determines a unique map \(P\to Q\) and a unique natural transformation \(\eta_0:\mathcal P'\to\mathcal Q\). The naturality of \(\mathcal P'f\) forces the image of \(P\to Q\) to lie in \(P_u\), because \(P_u\) is defined by the consistency equations imposed in Step 3. Hence \(P\to Q\) factors uniquely through \(P_u\), giving a map \(P\to P_u\) and a natural transformation \(\eta:\mathcal P'\to\mathcal P_u\) fitting into pullback squares
\[\begin{tikzcd}\mathcal P'(x)\dar\ar[bend left]{rr}{\eta_{0}}\rar{\eta}&\mathcal P_{u}(x)\rar\dar &\mathcal Q(x)\dar
\\ P\rar &P_{u}\rar &Q\end{tikzcd}\]

The identities defining \(P_u\) imply that
\[
        \mathcal P_u f\circ \eta=\mathcal P'f .
\]
Thus \(\eta:\mathcal P'\to\mathcal P_u\) is a map of extended tropical structures, and it is compatible with the chosen morphisms from \(\mathcal P_f\). The uniqueness of \(\eta_0\) in Claim~\ref{unique inclusion} implies uniqueness of \(\eta\). Therefore \(\iota_u:\mathcal P_f\to\mathcal P_u\) is the universal extension of the tropical structure of \(f\).

\end{proof}

\begin{remark}\label{pu char}
The universal extension \(\iota_u:\mathcal P_f\longrightarrow\mathcal P_u\) has a concrete  description.  Its base polytope \(P_u\) parametrizes the tropical deformations of the tropical curve \(\totb f\).  The chosen morphism \(\iota_u\) determines a distinguished point of \(P_u\), corresponding to the original tropical curve.

For each point \(x\) in a smooth stratum of \(\totl{\ex C}\), the identification \(\mathcal P_u(x)=P_u\) gives a map
\[
        A_x:P_u\longrightarrow \mathcal P(f(x)),
\]
equal to \(\mathcal P_u f:\mathcal P_u(x)\longrightarrow\mathcal P(f(x))\) under this identification.  The map \(A_x\) records the possible tropical position of the image of \(x\).  For each node \(e\) of \(\totl{\ex C}\), the construction gives a map
\[
        \rho_e:P_u\longrightarrow I_e:=[0,\infty),
\]
recording the possible length of the internal edge of \(\totb f\) corresponding to \(e\).  The local tropical polytope over the node is then the pullback
\[
\begin{tikzcd}
\mathcal P_u(e)\rar\dar & {[}0,\infty)^2\dar{a+b} \\
P_u\rar{\rho_e} & {[}0,\infty).
\end{tikzcd}
\]

Suppose that the nodes of \(\totl{\ex C}\) are \(e_1,\dotsc,e_n\), and choose points \(x_1,\dotsc,x_m\) in the smooth strata of \(\totl{\ex C}\), with at least one \(x_i\) in each smooth stratum.  Then the map
\[
        (\rho_{e_1},\dotsc,\rho_{e_n},A_{x_1},\dotsc,A_{x_m}):
        P_u\longrightarrow
        \prod_{j=1}^n I_{e_j}\times\prod_{i=1}^m\mathcal P(f(x_i))
\]
is injective and identifies \(P_u\) with a subpolytope of the target.  This subpolytope is cut out by integral-affine consistency equations: roughly, the equations assert that the proposed vertex positions and edge lengths actually fit together to define a tropical curve compatible with the tropical part of \(f\).

In the special case that \(f\) is a curve in a smooth manifold, each polytope \(\mathcal P(f(x_i))\) is a point.  Hence the edge-length map
\[
        (\rho_{e_1},\dotsc,\rho_{e_n}):P_u\longrightarrow[0,\infty)^n
\]
is an integral-affine isomorphism.  Thus, in the smooth case, the universal tropical deformation polytope is just the usual orthant of node-smoothing parameters.
\end{remark}

\section{Families of curves with universal tropical structure}

We now realize the universal extension constructed in Theorem~\ref{mts} by an actual family of curves. The first lemma constructs a family whose tropical structure is universal at the given curve \(f\). We then show that, after restricting the family if necessary, this property is open in the family. Finally, we prove the local mapping property used later in the construction of core families.

 \begin{defn}\label{uts def}
We shall say that a family \(\hat f\) has universal tropical structure at a curve \(f\) if the extension \(\mathcal P_f\to\mathcal P'\) determined by the inclusion of \(f\) as a fibre of \(\hat f\) is isomorphic to the universal extension \(\mathcal P_f\to\mathcal P_u\). A family has universal tropical structure if it has universal tropical structure at every curve it contains.
\end{defn}

\begin{lemma}\label{uts family}
Let \(f\) be a curve in \(\hat{\ex B}\) with domain \(\ex C\neq\ex T\). Then there exists a family of curves
\[
\begin{tikzcd}
\ex C(\hat f)\dar\rar{\hat f} & \hat{\ex B}\dar \\
\ex F(\hat f)\rar & \ex G
\end{tikzcd}
\]
containing \(f\), with universal tropical structure at $f$.  Moreover, the smooth part of the stratum of \(\ex F(\hat f)\) containing \(f\) is a single point.
\end{lemma}

\begin{proof}

The construction is the usual local construction of a family of curves, with the base polytope chosen to be \(P_u\). The edge-length maps \(\rho_e:P_u\to[0,\infty)\) determine the smoothing parameters at the nodes, while the maps \(A_x:P_u\to\mathcal P(f(x))\) determine the tropical part of the map to \(\hat{\ex B}\) on smooth components. We construct the domain first, then define the map \(\hat f\), and finally check that the induced extension of tropical structure is the universal one.

  Let \(p_f\in P_u\) be the distinguished point determined by the  morphism \(\iota_u:\mathcal P_f\to\mathcal P_u\). The base of the family will be an open subset of \(\et m{P_u}\). Choose a point \(\star_f\in \et m{P_u}\) with tropical part \(p_f\); the fibre over \(\star_f\) will be the original curve \(f\).

\smallskip
\noindent\emph{The domain.}

 Choose a finite collection of coordinate charts $U_{i}$ on the domain $\ex C$ of $f$ so that $U_{i}$ with its complex structure is either 
\begin{enumerate}
\item\label{t1} An open subset of $\et 1{[0,l]}$ in the form 
of 
\[1>\abs{\tilde z}>r\e l\]
\item \label{t2} an open subset of $\et 1{[0,\infty)}$ in the form 
\[1>\abs{\tilde z}\]
\item\label{t3} or an open subset of $\mathbb C$,
\end{enumerate}
so that
\begin{itemize}\item any nonempty  intersection of two different charts involves at least one chart of type \ref{t3}, 
\item and each  chart of type \ref{t3}  which intersects a chart of type \ref{t1} does so either  on the subset
\[1>\abs{\tilde z}>\frac 12\]
or on the subset
\[2r\e l>\abs{\tilde z}> r\e l\]
but does not intersect both of these subsets.

\end{itemize}

\

We shall extend these coordinate charts $U_{i}$ to coordinate charts $\tilde U_{i}$ on a family $\hat f$ of curves as follows:

In the case that $U_{i}$ is a chart of type \ref{t1}, proceed as follows: 
Denote by $V$ the open subset of $\et 2{[0,\infty)^{2}}$ in the form of 
\[V:=\left\{(\tilde z_{1},\tilde z_{2}): \ \abs{\tilde z_{1}}<1,\  \abs{\tilde z_{2}}<1,\   \abs{\tilde z_{1}\tilde z_{2}}<\frac 14\right\}\subset \et 2{[0,\infty)^{2}}\]
There is a submersion 
\[V\longrightarrow\et 1{[0,\infty)}\]
in the form of $(\tilde z_{1},\tilde z_{2})\mapsto \tilde z_{1}\tilde z_{2}$.

Suppose that $U_{i}$ is the coordinate chart containing the node $e$. Remark \ref{pu char} implies that  for each  node $e$ of $\totl{\ex C}$ there is a pullback diagram

\[\begin{tikzcd}\dar \mathcal P_{u}(e)\rar &{}[0,\infty)^{2}
\dar{ a+b}
\\P_{u}\rar{\rho_{e}} &{}[0,\infty)\end{tikzcd}\] 

There is  a  monomial map 
\[q_{i}:\et m{P_{u}}\longrightarrow \et 1{[0,\infty)}\]
 with tropical part equal to $\rho_{e}$ so that the fiber of $V$ over the image of $\star_f$ is isomorphic to $U_{i}$ (with its complex structure). Define $\ex F(\hat f)\subset \et m{P_{u}}$ to be the subset where $\abs {q_{i}}<\frac 14$  for all such maps $q_{i}$. Then define $\tilde U_{i}$ via  the pullback diagram 
\[\begin{tikzcd} \tilde U_i \rar\dar & V\dar
\\ \ex F(\hat f)\rar{q_i} & \et 1{[0,\infty)}\end{tikzcd}\]
and give $\tilde U_{i}\longrightarrow \ex F(\hat f)$ the fiberwise complex structure pulled back from $V\longrightarrow \et 1{[0,\infty)}$. Note that $U_{i}$ is the fiber of $\tilde U_{i}$ over $\star_f\in \ex F(\hat f)$.

\smallskip

If $U_{i}$ is any other coordinate chart (of the form \ref{t2} or \ref{t3}), then define $\tilde U_{i}$ as $U_{i}\times \ex F(\hat f)$ with the obvious projection to $\ex F(\hat f)$ and the fiberwise complex structure from $U_{i}$. Similarly, if $U_{I}$ indicates the intersection of $U_{i}$ for all $i\in I$, and $\abs{I}>1$, then define $\tilde U_{I}:=U_{I}\times \ex F(\hat f)$ with the obvious projection to $\ex F(\hat f)$.

\smallskip

To define a family of curves $\ex C(\hat f)\longrightarrow \ex F(\hat f)$ using the system of coordinate charts $\tilde U_{I}$, we still need to define transition maps. If $I\supset I'$, then there is a transition map 
\[\psi_{I,I'}:U_{I}\longrightarrow U_{I'}\]
If $U_{I'}$ is not a chart of type \ref{t1}, then there is an obvious induced map
\[\tilde\psi_{I,I'}:\tilde U_{I}:=U_{I}\times\ex F(\hat f)\xrightarrow{(\id,\psi_{I,I'})}U_{I'}\times \ex F(\hat f):= \tilde U_{I'}\]
which is a fiberwise holomorphic isomorphism onto its image. If on the other hand $U_{i}$ is a chart of type \ref{t1} and $I$ contains $i$, and $\psi_{I,i}$ has image within $1>\abs{\tilde z}>\frac 12$, then define $\tilde\psi_{I,i}$
as the unique map  
\[\tilde \psi_{I,i}:\tilde U_{I}\longrightarrow \tilde U_{i}\]
so that the following two diagrams commute.
\[\begin{tikzcd}[column sep=tiny]\tilde U_{I}\ar{dr}\ar{rr}{\tilde\psi_{I,i}}&&\ar{dl}\tilde U_{i}
\\ &\ex F(\hat f)\end{tikzcd}\]
\[\begin{tikzcd}\dar\tilde U_{I}\rar{\tilde \psi_{I,i}}&\tilde U_{i}\rar& V\dar{\tilde z_{1}}
\\ U_{I}\rar{\psi_{I,i}}&U_{i}\rar{\tilde z}&\et 11
\end{tikzcd}\]
 Similarly, if $\psi_{I,i}$ has image with $ 2r\e l>\abs{\tilde z}>r\e l$, then define $\tilde \psi_{I,i}$ to be the unique map so that the following diagram and the first of the above diagrams commute. 
\[\begin{tikzcd}\dar\tilde U_{I}\rar{\tilde \psi_{I,i}}&\tilde U_{i}\rar& V\dar{\tilde z_{2}}
\\ U_{I}\rar{\psi_{I,i}}&U_{i}\rar{r\e l\tilde z^{-1}}&\et 11
\end{tikzcd}\]
 
 As in each case, these coordinates of $V$ are (fiberwise) holomorphic functions on $\tilde U_{i}$, these transition maps are fiberwise holomorphic. $\tilde \psi_{I,i}$ is also an isomorphism onto its image. 
 
 With these definitions, compatibility of the coordinate maps $\tilde \psi_{I',I}$ follows from compatibility of the coordinate maps $\psi_{I',I}$. As the map $\prod_{i\in I}\psi_{I,i}$ is proper, it follows that $\prod_{i\in I}\tilde \psi_{I,i}$ is proper too. Therefore $\tilde U_{I}$ with these transition maps define a coordinate system on some family of curves $\ex C(\hat f)\longrightarrow \ex F(\hat f)$.  As the fiber of $\tilde U_{i}$ over $\star_f\in\ex F(\hat f)$ is always $U_{i}$ and the restriction of $\tilde \psi_{I',I}$ to this fiber is $\psi_{I',I}$, the fiber of $\ex C(\hat f)$ over $\star_f\in\ex F(\hat f)$ is isomorphic to $\ex C$.
 
 \begin{remark}\label{equivariant coordinates}
 In the proof of Theorem~\ref{G uts family} below, we shall use that $\{\tilde U_{I}\}$ is an equivariant set of coordinate charts on $\ex C(\hat f)\longrightarrow \ex F(\hat f)$ in a sense defined precisely in \cite[Appendix A]{cem}. In particular, the projections to $\et m{P_{u}}$ and the maps $\tilde\psi_{I,I'}$ are compatible with multiplication of exploded coordinates in the source and target by elements of $\mathbb C^{*}\e{\mathbb R}$. 
\end{remark}

\smallskip
\noindent\emph{The map to \(\hat{\ex B}\).}
Now we shall define a map \(\hat f:\ex C(\hat f)\to\hat{\ex B}\) so that
\begin{enumerate}
\item\label{fc} \(\hat f\) restricted to \(\ex C\subset\ex C(\hat f)\) is equal to \(f\);
\item\label{tropc} the tropical structure of \(\hat f\) restricted to \(\ex C_T\) is the universal extension, namely
\[
\begin{array}{ccc}
(\ex C_T,\mathcal P_u)&\xrightarrow{\mathcal P_u f}&(\hat{\ex B}_T,\mathcal P) \\
\downarrow&&\downarrow \\
(*,P_u)&\longrightarrow&(\ex G_T,\mathcal P).
\end{array}
\]
\end{enumerate}
 
We have constructed $\hat f$ so that the tropical structure of $\ex C(\hat f)\longrightarrow \ex F(\hat f)$ restricted to $\ex C_{T}$ is equal to the left hand side of the above diagram. The above condition (\ref{tropc}) on the tropical structure of $\hat f$ translates to requiring that restricted to any coordinate chart on $\ex C(\hat f)$ containing a point $x\in\totl{\ex C}$ in its interior stratum, the tropical part of $\hat f$ in this coordinate chart is equal to $\mathcal P_{u}f$. The fact that $\mathcal P_{u}f$ is a natural transformation implies that condition (\ref{tropc}) does not depend on the point $x\in\totl{\ex C}$ chosen, and that the restriction of a map satisfying condition (\ref{tropc}) on one coordinate chart to another coordinate chart will still satisfy condition (\ref{tropc}). 

Choose a map $\et m{ P_{u}}\longrightarrow \ex G$ with the correct tropical part which sends $\star_f\in\ex F(\hat f)$ to the point in  $\ex G$ whose fiber in $\hat{\ex B}$ contains the image of $f$.

As $\mathcal P_{u}$ extends the tropical structure of $f$, the tropical condition (\ref{tropc}) is compatible with the condition \ref{fc} that $\hat f$ restricted to $\ex C$ is equal to $f$. Therefore, around any point in $\ex C$, there locally exists a map obeying the above two  conditions and projecting to the given map $\et m{ P_{u}}\longrightarrow \ex G$. Any two such maps differ by exponentiation of a section of the pullback of $T_{vert}\hat{\ex B}$ which vanishes on $\ex C\subset \ex C(\hat f)$. Exponentiating  such a vector field always preserves conditions (\ref{fc}) and (\ref{tropc}). Therefore, we may extend a map satisfying conditions (\ref{fc}) and (\ref{tropc}) on some set of coordinate charts to the next coordinate chart by using a cutoff function to interpolate between the old map and a new map on the next coordinate chart. 

In this way, we may construct a map 
\[\begin{array}{ccc}\ex C(\hat f)&\xrightarrow{\hat f} &\hat {\ex B}
\\\downarrow&&\downarrow
\\\ex F(\hat f)&\longrightarrow&\ex G\end{array}\]
 extending $f$ with tropical structure which restricts to $\ex C_{T}$ to be the universal extension of the tropical structure of $f$.
 
 \end{proof}

\begin{lemma}\label{open uts}Suppose that $\hat f$ is a family of curves with universal tropical structure at some curve $f$. Then $\hat f$ also has universal tropical structure at any curve $f'$ in a stratum of $\ex F(\hat f)$ with closure containing $f$. \end{lemma}

\begin{proof}

\

\smallskip
\noindent\emph{Setup and the map \(P'\to P_u^{f'}\).}

Let \(\mathcal R\) denote the extension of the tropical structure of \(f'\) obtained by restricting the tropical structure of \(\hat f\) to \(f'\), and let \(P'\) denote its base polytope. Let \(\mathcal P_u^{f'}\) denote the universal extension of the tropical structure of \(f'\), with base polytope \(P_u^{f'}\). Similarly, let \(\mathcal P_u^f\) denote the universal extension of the tropical structure of \(f\), with base polytope \(P_u^f\). Since \(\hat f\) has universal tropical structure at \(f\), we identify the extension determined by \(\hat f\) at \(f\) with \(\mathcal P_u^f\).

The universal property of \(\mathcal P_u^{f'}\) gives a map of extended tropical structures \(\mathcal R\to\mathcal P_u^{f'}\), and hence a map of base polytopes \(P'\to P_u^{f'}\). It suffices to prove that this map is an isomorphism.

The idea of the proof is as follows. Parallel transport from \(f'\) to \(f\) identifies \(P'\) with a face of \(P_u^f\). Similarly, the unconstrained parameter polytope \(Q'\) used to construct the universal extension of \(f'\) embeds as a face of the corresponding polytope \(Q\) for \(f\). The universal polytope \(P_u^{f'}\subset Q'\) maps to \(Q\), and we must show that its image is exactly the face \(P'\subset P_u^f\). This is immediate when \(f'\) lies in the same stratum as \(f\). In general, we reduce to the codimension-one case and check that \(P_u^{f'}\), viewed inside \(Q\), satisfies the integral-affine equations defining \(P_u^f\). These equations are of two kinds: loop equations on smooth components and compatibility equations at nodes.

\smallskip
\noindent\emph{Comparison with the universal data for \(f\).}

The curve \(f'\) is in a stratum of \(\ex F(\hat f)\) with closure containing \(f\) if and only if there exists a path \(\gamma\) joining \(f'\) to \(f\) in \(\totl{\ex F(\hat f)}\) along which parallel transport of tropical structure is defined. Parallel transport along \(\gamma\) gives an inclusion of \(P'\) as a face of \(P_u^f\). Similarly, parallel transport along a lift of this path joining \(x'\in\totl{\ex C(f')}\) to \(x\in\totl{\ex C(f)}\) gives an inclusion of \(\mathcal R(x')\) as the face of \(\mathcal P_u^f(x)\) which is the inverse image of \(P'\subset P_u^f\).

We shall use the notation established in the proof of Theorem~\ref{mts}. Let \(Q\) be the ambient parameter polytope used there for the curve \(f\). The path \(\gamma\) may be lifted to paths \(\gamma_i\) in \(\totl{\ex C(\hat f)}\), along which parallel transport is defined, joining points \(p_i'\in\totl{\ex C(f')}\) to the chosen points \(p_i\in\totl{\ex C(f)}\). Parallel transport along \(\gamma\), \(\gamma_i\), and \(\hat f_T\gamma_i\) fits into the following commutative diagram.
\begin{equation}\label{cd1}
\begin{tikzcd}[column sep=large]
\mathcal P(f'(p_i'))\rar{\mathcal P(\hat f_T\gamma_i)}&\mathcal P(f(p_i)) \\
\mathcal R(p_i')\dar\uar{\mathcal P\hat f}\rar{\mathcal P(\gamma_i)}&\mathcal P_u^f(p_i)\uar{\mathcal P\hat f}\dar \\
P'\rar{\mathcal P(\gamma)}&P_u^f
\end{tikzcd}
\end{equation}

Define the polytope \(Q'\) to be the product of \(\mathcal P(f'(p_i'))\) for each \(i\), together with one factor \(I_{e'}:=[0,\infty)\) for each node \(e'\) of \(\totl{\ex C(f')}\):
\[
        Q':=\prod_i\mathcal P(f'(p_i'))\times\prod_{e'}I_{e'} .
\]
This is the analogue for \(f'\) of the polytope \(Q\) defined in equation~(\ref{Qdef}) for \(f\). Every node of \(\totl{\ex C(f')}\) may be followed along \(\gamma\) to identify it with a node of \(\totl{\ex C(f)}\). Thus the nodes of \(\totl{\ex C(f')}\) identify with a subset of the nodes of \(\totl{\ex C(f)}\). There is therefore an inclusion of \(Q'\) as a face of \(Q\), induced by the maps \(\mathcal P(\hat f_T\gamma_i)\) on the vertex-position factors, the identity maps on the \(I_e\) factors corresponding to identified nodes, and the zero map on the \(I_e\) factors of \(Q\) corresponding to nodes of \(\totl{\ex C(f)}\) which do not come from nodes of \(\totl{\ex C(f')}\).

There is then a commutative diagram
\begin{equation}\label{cd2}
\begin{array}{ccc}
Q'&\longrightarrow &Q \\
\uparrow &&\uparrow \\
P'&\longrightarrow &P_u^f .
\end{array}
\end{equation}
Here the map \(P'\longrightarrow Q'\) is defined analogously to the map \(P\longrightarrow Q\) from Claim~\ref{unique inclusion} in the proof of Theorem~\ref{mts}.

We may choose the paths \(\gamma_i\) so that if \(p_i'\) and \(p_j'\) are in the same smooth component of \(\ex C(f')\), then \(p_i'=p_j'\). Restrict to the subpolytope of \(Q'\) defined by requiring the coordinates corresponding to such repeated factors \(\mathcal P(f'(p_i'))\) and \(\mathcal P(f'(p_j'))\) to agree, and then repeat the construction from the proof of Theorem~\ref{mts}. This constructs \(P_u^{f'}\) as a subpolytope of \(Q'\). The map \(P'\longrightarrow Q'\) has image in \(P_u^{f'}\subset Q'\).

The lower and right-hand arrows in diagram~(\ref{cd2}) are inclusions of subpolytopes whose images are defined by integral-affine equations. It follows that the map \(P'\longrightarrow P_u^{f'}\subset Q'\) is an inclusion of \(P'\) as a subpolytope of \(P_u^{f'}\), again with image defined by integral-affine equations. To prove that \(P'\longrightarrow P_u^{f'}\) is an isomorphism, it therefore suffices to prove that this map is surjective. Equivalently, it suffices to prove that the image of \(P_u^{f'}\subset Q'\) in \(Q\) under the top arrow of diagram~(\ref{cd2}) is contained in the image of \(P'\) in \(P_u^f\subset Q\).

\smallskip
\noindent\emph{Reduction to a codimension-one face.}

If \(f'\) is in the same stratum as \(f\), then the constructions of \(P_u^{f'}\) and \(P_u^f\) use isomorphic polytopes \(Q'\) and \(Q\), and isomorphic integral-affine relations to define the subpolytopes \(P_u^{f'}\) and \(P_u^f\). In this case \(P_u^{f'}\), \(P_u^f\), and \(P'\) have the same dimension, so the inclusion \(P'\longrightarrow P_u^{f'}\) is an isomorphism.

We may therefore reduce to the case that \(P'\) is one dimension smaller than \(P_u^f\). Since any face of \(P_u^f\) of smaller dimension may be reached by taking a sequence of codimension-one faces, the general case follows by repeating the argument below. We must prove that the image of \(P_u^{f'}\) in \(Q\) is contained in the image of \(P'\) in \(Q\).

Because \(P'\) is one dimension smaller than \(P_u^f\), the curve \(f'\) is in a different stratum from \(f\). Thus the image of the curve \(f\) in \(P_u^f\subset Q\) cannot be contained in \(Q'\subset Q\): this would require either some \(p_i\) to map outside the interior of \(\mathcal P(f(p_i))\), or some edge of \(\totb{\ex C(f)}\) to have zero length. Therefore \(Q'\cap P_u^f\) is at least one dimension smaller than \(P_u^f\). Since the image of \(P'\) in \(Q\) is one dimension smaller than \(P_u^f\) and is contained in \(Q'\cap P_u^f\), the image of \(P'\) in \(Q\) is equal to \(Q'\cap P_u^f\). Here we use that all these polytopes are obtained by intersecting \(Q\) with affine subspaces.

It remains to show that the image of \(P_u^{f'}\) in \(Q\) is contained in \(P'=Q'\cap P_u^f\). Equivalently, we must show that the image of \(P_u^{f'}\) in \(Q\) is contained in \(P_u^f\subset Q\). We shall do this by checking that \(P_u^{f'}\), viewed inside \(Q\), satisfies the integral-affine equations defining \(P_u^f\).

\smallskip
\noindent\emph{Loop equations.}

First consider the equations coming from loops in smooth components. Given a loop \(\gamma_0\) in a smooth component of \(\ex C(f)\), starting and ending at \(p_i\), there is a map
\[
        \mathcal P(f_T\gamma_0):\mathcal P(f(p_i))\longrightarrow\mathcal P(f(p_i)).
\]
The polytope \(P_u^f\) is contained in the subset of \(Q\) whose corresponding coordinates are fixed by \(\mathcal P(f_T\gamma_0)\).

There exists a loop \(\gamma_0'\) in a smooth component of \(\ex C(f')\), starting and ending at \(p_i'\), such that \(\gamma_0'\) followed by \(\gamma_i\) is homotopic to \(\gamma_i\) followed by \(\gamma_0\). Therefore
\[
        \mathcal P(\hat f_T\gamma_i)\circ\mathcal P(\hat f_T\gamma_0')
        =
        \mathcal P(\hat f_T\gamma_0)\circ\mathcal P(\hat f_T\gamma_i).
\]
Since \(P_u^{f'}\) is contained in the subset of \(Q'\) whose corresponding coordinates are fixed by \(\mathcal P(\hat f_T\gamma_0')\), and since the inclusion \(Q'\longrightarrow Q\) on the relevant factor is \(\mathcal P(\hat f_T\gamma_i)\), the image of \(P_u^{f'}\) in \(Q\) has coordinates fixed by \(\mathcal P(f_T\gamma_0)\). Hence the image of \(P_u^{f'}\) in \(Q\) is contained in the polytope \(Q_1\subset Q\) from the proof of Theorem~\ref{mts}.

\smallskip
\noindent\emph{Node equations.}

It remains to check the node equations defining \(P_u^f\subset Q_1\). These equations come from pairs of paths \(h_1\) and \(h_2\) joining \(p_i\) and \(p_j\), respectively, to a node \(e\) of \(\totl{\ex C(f)}\). These paths were denoted by \(\gamma\) and \(\gamma'\) in the proof of Theorem~\ref{mts}; we use different notation here to avoid conflict with the path \(\gamma\) in the base.

Lift \(\gamma\) to a path \(\gamma_e\) in \(\totl{\ex C(\hat f)}\) ending at \(e\), and along which parallel transport is defined. Let \(e'\in\totl{\ex C(f')}\) be the starting point of \(\gamma_e\). There are paths \(h_1'\) and \(h_2'\) in \(\totl{\ex C(f')}\), joining \(p_i'\) and \(p_j'\) to \(e'\), such that \(h_1'\) followed by \(\gamma_e\) is homotopic to \(\gamma_i\) followed by \(h_1\), and \(h_2'\) followed by \(\gamma_e\) is homotopic to \(\gamma_j\) followed by \(h_2\). Parallel transport along the images of these paths under \(\totl{\hat f}\) gives the following commutative diagram:
\begin{equation}\label{ocd1}
\begin{tikzcd}[column sep=large]
\mathcal P(f'(p_i'))\dar[swap]{\mathcal P(\hat f_T h_1')}\rar{\mathcal P(\hat f_T\gamma_i)}&\mathcal P(f(p_i))\dar{\mathcal P(\hat f_T h_1)} \\
\mathcal P(f'(e'))\rar{\mathcal P(\hat f_T\gamma_e)}&\mathcal P(f(e)) \\
\mathcal P(f'(p_j'))\rar{\mathcal P(\hat f_T\gamma_j)}\uar{\mathcal P(\hat f_T h_2')}&\mathcal P(f(p_j))\uar[swap]{\mathcal P(\hat f_T h_2)}
\end{tikzcd}
\end{equation}

There is a unique map \(P_u^{f'}\longrightarrow I_e=[0,\infty)\) such that the local polytope \(\mathcal P_u^{f'}(e')\) is obtained by the pullback diagram
\[
\begin{tikzcd}
\mathcal P_u^{f'}(e')\rar\dar
    & {[}0,\infty)^2\dar{a+b} \\
P_u^{f'}\rar
    & {[}0,\infty).
\end{tikzcd}\]
If \(e'\) is not a node of \(\totl{\ex C(f')}\), then this map has image \(0\) in \([0,\infty)\), and \(\mathcal P_u^{f'}(e')\) is isomorphic to \(P_u^{f'}\). Since \(\mathcal Q(e)\) is constructed similarly as a pullback of \([0,\infty)^2\longrightarrow [0,\infty)\) over the corresponding coordinate of \(Q\), we get a canonical map \(\mathcal P_u^{f'}(e')\longrightarrow\mathcal Q(e)\) so that the following diagram commutes:
\begin{equation}\label{ocd2}
\begin{tikzcd}
\mathcal P_u^{f'}(p_i')=P_u^{f'}\rar\dar[swap]{\mathcal P_u^{f'}(h_1')} &Q=\mathcal Q(p_i)\dar{\mathcal Q(h_1)} \\
\mathcal P_u^{f'}(e')\rar &\mathcal Q(e) \\
\mathcal P_u^{f'}(p_j')=P_u^{f'}\rar\uar{\mathcal P_u^{f'}(h_2')} &Q=\mathcal Q(p_j)\uar[swap]{\mathcal Q(h_2)}
\end{tikzcd}
\end{equation}

Define \(A_{h_i}\) analogously to the definition of \(A_{\gamma}\) given in diagram~(\ref{Adef}) on page~\pageref{Adef}:
\begin{equation}\label{ocd3}
\begin{tikzcd}
\mathcal P(e)\rar{\eta_0}\ar[bend left]{rr}{\mathcal Pf}&\mathcal Q(e)\rar{A_{h_i}}&\mathbb R^n\supset\mathcal P(f(e)) \\
&\mathcal Q(h_i(0))\uar{\mathcal Q(h_i)}\rar{\mathcal Qf}&\mathcal P(f(h_i(0)))\uar[swap]{\mathcal P(f_T(h_i))}
\end{tikzcd}
\end{equation}

Consider the following diagram:
\[
\begin{tikzcd}
\ar[bend left]{rrr}{\mathcal P(\hat f_T\gamma_i)}
\mathcal P(f'(p_i'))\dar{\mathcal P(\hat f_T h_1')}
    &\lar{\mathcal P_u^{f'}f'}\mathcal P_u^{f'}(p_i')=P_u^{f'}\rar\dar[swap]{\mathcal P_u^{f'}(h_1')}
    &Q=\mathcal Q(p_i)\dar{\mathcal Q(h_1)}\rar{\mathcal Qf}
    &\mathcal P(f(p_i))\dar{\mathcal P(\hat f_T h_1)} \\
\ar[bend right]{rrr}{\mathcal P(\hat f_T\gamma_e)}
\mathcal P(f'(e'))
    &\lar{\mathcal P_u^{f'}f'}\mathcal P_u^{f'}(e')\rar
    &\mathcal Q(e)\rar{A_{h_1}}
    &\mathbb R^n\supset\mathcal P(f(e)).
\end{tikzcd}
\]
Diagram~(\ref{ocd1}) implies that the outer loop commutes. The naturality of \(\mathcal P_u^{f'}f'\) implies that the left-hand square commutes. The middle square commutes by diagram~(\ref{ocd2}), the right-hand square commutes by diagram~(\ref{ocd3}), and the top loop commutes because the definition of the map \(P_u^{f'}\subset Q'\longrightarrow Q\) uses the map \(\mathcal P(\hat f_T\gamma_i)\) on the relevant factor projected onto by \(\mathcal P_u^{f'}f'\) and \(\mathcal Qf\). Therefore the bottom loop commutes when \(\mathcal P_u^{f'}(e')\) is restricted to the image of \(\mathcal P_u^{f'}(p_i')\). The bottom loop also commutes when \(\mathcal P_u^{f'}(e')\) is restricted to the complementary image of the original tropical polytope \(\mathcal P(e')\). Hence the bottom loop of the diagram above commutes.

Repeating the same argument for \(h_2\) gives that the following diagram commutes:
\[
\begin{tikzcd}[column sep=large]
&\mathcal Q(e)\ar[end anchor=40]{dr}{A_{h_1}} \\
\mathcal P_u^{f'}(e')\ar{dr}\ar{ur}\rar{\mathcal P_u^{f'}f'}&\mathcal P(f'(e'))\rar{\mathcal P(\hat f_T\gamma_e)}&\mathcal P(f(e))\subset\mathbb R^n \\
&\mathcal Q(e)\ar[end anchor=-40]{ur}[swap]{A_{h_2}}
\end{tikzcd}
\]

It follows that the image of \(\mathcal P_u^{f'}(e')\) in \(\mathcal Q(e)\) is contained in the subset where \(A_{h_1}=A_{h_2}\). Since \(P_u^f\) is the subset of \(Q_1\) over which \(A_{h_1}=A_{h_2}\) for all such pairs of paths, the image of \(P_u^{f'}\) in \(Q\) is contained in \(P_u^f\). As noted above, this implies that \(P'\to P_u^{f'}\) is an isomorphism. Therefore \(\hat f\) has universal tropical structure at \(f'\).

\end{proof}

\begin{lemma}\label{uts map} Let $\hat f$ be a family of curves with universal tropical structure at some curve $f$.  Given any other family of curves $\hat h$ containing $f$, by restricting to a neighborhood of $f$ in $\hat h$ there exists a map 
\[\begin{array}{ccc}\ex C(\hat h)&\xrightarrow{\Phi} &\ex C(\hat f)
\\\downarrow&&\downarrow
\\\ex F(\hat h)&\longrightarrow &\ex F(\hat f)
 \end{array}\]
which is the identity map on $\ex C(f)$ within $\ex C(\hat h)$ and $\ex C(\hat f)$, and so that in a metric on $\hat{\ex B}$, the distance between the maps $\hat h$ and $\hat f\circ \Phi$ is bounded.

\end{lemma}

\begin{proof}

The bounded-distance condition forces \(\hat h\) and \(\hat f\circ\Phi\) to have the same tropical part. Thus the tropical part of such a \(\Phi\) must be the unique map from the extension determined by \(\hat h\) to the universal extension determined by \(\hat f\). We use this unique map of extended tropical structures to prescribe the tropical part of \(\Phi\).

In particular, given any point $x$ in $\totl{\ex C(f)}$, we get a specification of what the tropical part of $\Phi$ should be restricted to coordinate charts on $\ex C(\hat h)$ and $\ex C(\hat f)$ containing $x$. The naturality of the map of extended tropical structures implies the following: given a map defined on one coordinate chart with the correct tropical part,   the restriction of this map to another coordinate chart  will again have the correct tropical part.   

Choose a map $\ex F(\hat h)\longrightarrow\ex F(\hat f)$ with the correct tropical part and sending $f$ to $f$. (If necessary, restrict to a neighborhood of $f$ in $\ex F(\hat h)$ so that such a map exists.)
Now construct a $\Phi$ so that: 
\begin{itemize}
\item the following diagram commutes
\[\begin{array}{ccc}\ex C(\hat h)&\xrightarrow{\Phi} &\ex C(\hat f)
\\\downarrow&&\downarrow
\\\ex F(\hat h)&\longrightarrow &\ex F(\hat f)
 \end{array}\]
\item $\Phi$ is the identity restricted to $\ex C(f)$; and
\item the tropical part of $\Phi$ in a coordinate chart containing $x\in \ex C(f)$ is the same as the map coming from the unique map of extended tropical structures.
\end{itemize}

Around $x$ in $\ex C(f)$, there locally  exists a map satisfying the above properties.  Any two such maps differ by the flow of a  vertical vector field which vanishes on $\ex C( f)$, and the flow of any such vertical vector field preserves the above properties.  If $\Phi$ is defined on some collection of coordinate charts on $\ex C(\hat f)$ intersecting $\ex C(f)$,  we can extend the domain of  definition of $\Phi$ to include the next coordinate chart by using a cutoff function to interpolate between the previously defined $\Phi$ and a map defined on the new coordinate chart. This constructs such a $\Phi$ on a neighborhood of $\ex C(f)$.

\end{proof}

\

We shall also need the following variant of Lemma~\ref{uts map}. It applies when the curve in the second family is a refinement of \(f\), or is obtained from \(f\) by adding external edges or bubble components. This variant is used in~\cite{evc} to define evaluation maps to Deligne--Mumford space, and in proving invariance of exploded Gromov--Witten invariants under refinement.

\begin{lemma}\label{uts map cor} 
 Let $\hat f$ be a family of curves with universal tropical structure at $f$.  Let $\hat h$ be a family of curves  containing  a  curve  $f'$ with a degree $1$ holomorphic map 
 \[\phi:\ex C( f')\longrightarrow \ex C(f) \] 
 so that $f'=f\circ\phi$. Then by restricting to a neighborhood of $f'$ in $\hat h$ there exists an extension of $\phi$ to a map 
\[\begin{array}{ccc}\ex C(\hat h)&\xrightarrow{\Phi} &\ex C(\hat f)
\\\downarrow&&\downarrow
\\\ex F(\hat h)&\longrightarrow &\ex F(\hat f)
 \end{array}\]
 so that in a metric on $\hat{\ex B}$, the distance between the maps $\hat h$ and $\hat f\circ \Phi$ is bounded.

\end{lemma}

\begin{proof}

We shall first construct the tropical part of $\Phi$. Then the existence of a map $\Phi$ satisfying the required properties will follow as in the proof of Lemma \ref{uts map}.

\

Denote the tropical structure of $\hat h$ restricted to $\ex C(f')_{T}$ by $\mathcal P_{0}$.
To construct the tropical part of $\Phi$, we need to construct a map from $\mathcal P_{0}$ to the universal extension of the tropical structure of $f$. To achieve this, it suffices to construct an extension $\mathcal P_{1}$ of the tropical structure of $f$ with a map from  $\mathcal P_{0}$.
\[\begin{tikzcd}[row sep=tiny](\ex C(f')_{T},\mathcal P_{0}) \ar{dd}{(\phi_{T},\eta)} \ar[bend left=15]{rrd}{(f'_{T},\mathcal P_{0}f')}
\\&&(\ex B_{T},\mathcal P)
\\(\ex C_{T},\mathcal P_{1})\ar[bend right=15]{rru}[swap]{(f_{T},\mathcal P_{1}f)} \end{tikzcd}\]
Denote by $P$ the polytope $\mathcal P_{0}(\ex F(f'))$. We shall construct $\mathcal P_{1}$ so that $\mathcal P_{1}(\ex F(f))=P$ too. 

We shall define $\mathcal P_{1}$ and the natural transformation $\eta$ separately on the inverse image of each point in $\ex C(f)$, and then show these constructions may be glued together.  

\

 For any point $x$ in a smooth stratum  of $\ex C(f)$, set $\mathcal P_{1}(x)$ equal to $P$. For any point $y$ in $\ex C(f')_{T}$ sent to $x$ by $\phi_{T}$, define \[\eta:\mathcal P_{0}(y)\longrightarrow \mathcal P_{1}(x)\]
 to be the projection to $P$ coming from the tropical structure of the projection $\ex C(\hat h)\longrightarrow\ex F(\hat h)$. For any smooth stratum $C$ of $\ex C(f)$, the fact that $\phi$ is holomorphic and degree $1$ implies that there exists a unique smooth component $C_{0}$ of $\phi^{-1}(C)\subset \ex C(f')$ on which  $\phi$ is injective. The image $\phi(C_{0})\subset C$ is then dense. The inverse image of $x\in C$ is either a single point in $C_{0}$ or a connected and closed union of strata of $\ex C(f')$ attached to $C_{0}$ which map constantly into $x$. If $y$ is a node or marked point of $\totl{\ex C(f')}$, the fact that the stratum corresponding to $y$ maps to $x$ implies that the map $\mathcal P_{0}f':\mathcal P_{0}(y)\longrightarrow \mathcal P(f(x))$ factors through projection to $P$ followed by a map $P\longrightarrow\mathcal P(f(x))$. The same holds for all points $y$ in the inverse image of $x$ and the inverse image of $x$ is connected, therefore the maps $P\longrightarrow \mathcal P(f(x))$ must be the same for each $y$ in the inverse image of $x$. We may therefore  define 
 \[\mathcal P_{1}f:\mathcal P_{1}(x)\longrightarrow \mathcal P(f(x))\] 
so that for any path $\gamma$ in $\phi^{-1}(x)$ along which parallel transport is defined,  the diagram 
\[\begin{tikzcd}\mathcal P_{0}(\gamma(0))\dar{\eta} \ar{dr} \rar{\mathcal P_{0}(\gamma)}&\mathcal P_{0}(\gamma(1))\ar{d}{\mathcal P_{0}f'}\ar{dl}
\\ P=\mathcal P_{1}(x)\rar{\mathcal P_{1}f}&\mathcal P (f(x))\end{tikzcd}\]
commutes. For $x$ a point in a smooth component of $\ex C(f)$,  we have now constructed $\mathcal P_{1}$ for the point $x$, and constructed  $\eta$ on $\phi^{-1}(x)$.

\

Suppose that $e$ is a node or puncture of $\totl{\ex C(f)}$. We shall now construct $\mathcal P_{1}(e)$ and $\eta$ on $\phi^{-1}(e)$. Choose an integral-affine identification of the edge of $\totb{\ex C(f)}$   corresponding to $e$ as a sub interval of $[0,\infty)$ with closure containing $0$. The inverse image of $e$ in $\ex C(f')$ is a connected, closed union of strata.  There is a finite collection of strata  $e_{1},\dotsc,e_{n}$ of $\phi^{-1}(e)$ with non-constant image in $[0,\infty)$. The fact that $\phi$ is a degree $1$ map  implies that the image of these strata do not intersect, and each maps isomorphically to a subinterval of $[0,\infty)$. Order the $e_{i}$ so that $e_{i}$ has image in $[0,\infty)$ before the image of $e_{i+1}$. In the case that $e$ is a node, these strata correspond to nodes $e_{1},\dotsc, e_{n}$. In the case that $e$ is a puncture, $e_{n}$ is a puncture and the other $e_{i}$ are nodes. 

According to Definition \ref{ets} part \ref{ets2},  corresponding to each node $e_{i}$, there is a pullback diagram
\[\begin{tikzcd}\mathcal P_{0}(e_{i})\rar \dar &{[}0,\infty)^{2}\dar{a+b}
\\ P\rar{\rho_{i}}&{[}0,\infty)\end{tikzcd}\]

If $e$ is a node, define $\mathcal P_{1}(e)$ via the pullback diagram
\[\begin{tikzcd}\mathcal P_{1}(e)\rar \dar &{[}0,\infty)^{2} \dar{a+b}
\\ P\rar{\sum_{i=1}^{n}\rho_{i}}&{[}0,\infty)\end{tikzcd}\]

Let $p\in P$ be the point corresponding to the tropical part of the inclusion $\ex F(f')\longrightarrow \ex F(\hat h)$. The length of the edge  in $\totb{\ex C(f')}$ corresponding to $e_{i}$ is $\rho_{i}(p)$, so the length of the edge corresponding to $e$ in $\totb{\ex C(f)}$ is equal to $\sum_{i}\rho_{i}(p)$. In particular, we may choose an isomorphism of $\mathcal P(e)$ with the fiber of $\mathcal P_{1}(e)\longrightarrow P$ over $p$.

If on the other hand, $e$ is an external edge, define $\mathcal P_{1}(e):=P\times [0,\infty)$. In this case, the fiber of $\mathcal P_{1}(e)$ over $p\in P$ is uniquely isomorphic to $\mathcal P(e)$.

In either case, for any $y$ in the inverse image of $e$, there exists a unique integral-affine map 
\[\eta:\mathcal P_{0}(y)\longrightarrow \mathcal P_{1}(e)\]
so that the following diagram commutes.

\[\begin{tikzcd}\mathcal P(y)\dar \rar{\mathcal P\phi}&\mathcal P(e)\dar
\\\mathcal P_{0}(y)\rar{\eta} \dar&\mathcal P_{1}(e)\ar{dl}
\\ P\end{tikzcd}\]

Given any path $\gamma$ in the inverse image of $e$ joining $y'$ to $y$, the following diagram commutes
\[\begin{tikzcd}\mathcal P(y')\ar[bend left]{rr}{\mathcal P\phi}\rar{\mathcal P(\gamma)} \dar &\mathcal P(y)\dar \rar{\mathcal P\phi}&\mathcal P(e)\dar
\\\mathcal P_{0}(y')\ar[bend left, dotted]{rr}\rar{\mathcal P_{0}(\gamma)} \ar{dr}&\mathcal P_{0}(y)\rar{\eta} \dar &\mathcal P_{1}(e)\ar{dl}
\\& P\end{tikzcd}\]
Therefore, the uniqueness of $\eta$ implies that the composition of the arrows in the middle row is $\eta:\mathcal P_{0}(y')\longrightarrow\mathcal P_{1}(e)$. Therefore, $\eta$ defines a natural transformation from $\mathcal P_{0}$ restricted to $\phi^{-1}_{T}(e)$ to the constant functor with image $\mathcal P_{1}(e)$. 

\

So far,  for every object $x$ in $\ex C_{T}( f)$  we have the following commuting diagram of functors from $\phi^{-1}_{T}(x)$
\[\begin{tikzcd}\mathcal P\rar{\mathcal P\phi} \dar&\mathcal P\circ\phi_{T} \dar
\\\mathcal P_{0}\rar{\eta}&\mathcal P_{1}\circ\phi_{T}\end{tikzcd}\]
In fact, the $\mathcal P_{1}\circ \phi_{T}$ in the above diagram may be regarded as a pushout in the following sense:
\begin{claim}\label{Ppushout} Given any constant functor $\mathcal P'$ from $\phi^{-1}_{T}(x)$ to a polytope $P'$  and a commuting diagram of functors on $\phi^{-1}_{T}(x)$,
\[\begin{tikzcd}\dar \mathcal P\rar{\mathcal P\phi}&\mathcal P\circ\phi_{T}\dar
\\ \mathcal P_{0}\rar & \mathcal P'\end{tikzcd}\]
there exists a unique map  
\[\mathcal P_{1}(x)\longrightarrow P'\]
corresponding to a natural transformation $\mathcal P_{1}\circ\Phi_{T}\longrightarrow \mathcal P'$ so that the following diagram commutes
\[\begin{tikzcd}\mathcal P\rar{\mathcal P\phi} \dar&\mathcal P\circ\phi_{T} \dar \ar[bend left=15]{ddr}
\\\mathcal P_{0}\rar{\eta}\ar[bend right=15]{rrd}&\mathcal P_{1}\circ\phi_{T}\ar[dotted]{dr}{\exists !}
\\ && \mathcal P'\end{tikzcd}\]

\end{claim}

To prove Claim \ref{Ppushout}, consider the diagram
\[\begin{tikzcd}\mathcal P(y)\rar\dar &\mathcal P(x)\dar
\\ \mathcal P_{0}(y)\rar &\mathcal P_{1}(x)\end{tikzcd}\]
for any object $y$ in $\phi^{-1}_{T}(x)$. The associated diagram featuring the affine spaces generated by the affine polytopes is a pushout diagram of integral linear spaces. This may be seen by considering the following  three cases:  If the dimensions of $\mathcal P(y)$ and $\mathcal P(x)$ are equal, then the map $\mathcal P_{0}(y)\longrightarrow \mathcal P_{1}(x)$ generates an isomorphism of integral-affine spaces. If $y$ is a node or edge and $x$ is a point in a smooth component, then $\mathcal P(y)\longrightarrow \mathcal P_{0}(y)$ is injective, and the integral-affine map generated by $\mathcal P_{0}(y)\longrightarrow \mathcal P_{1}(x)$ is the quotient map which sends the span of  the image of $\mathcal P(y)$ to a single point.  If $y$ is a point in  a smooth component, and $x$ is a node or edge, then the integral-affine space spanned by $\mathcal P_{1}(x)$ is the span of the complementary images of $\mathcal P_{0}(y)$ and $\mathcal P(x)$. In each of the above three cases, the diagram of integral-affine spaces  generated by the above diagram is a pushout diagram.

 Therefore, if we embed $P'$ into $\mathbb R^{n}$, there exists a unique map $\mathcal P_{1}(x)\longrightarrow \mathbb R^{n}$ so that the following diagram commutes  
\[\begin{tikzcd}\mathcal P_{1}(x)\ar{dr}&\lar \mathcal P(x)\dar
\\ \mathcal P_{0}(y)\uar{\eta} \rar & \mathbb R^{n}\end{tikzcd}\]
The uniqueness of this map and the fact that $\phi^{-1}(x)$ is connected implies that the same map is obtained using any $y$ in the inverse image of $x$.
As the union of the images of $\eta$ is all of $\mathcal P_{1}(x)$, this map $\mathcal P_{1}(x)\longrightarrow \mathbb R^{n}$ must have image inside $P'$. This constructs the unique map $\mathcal P_{1}(x)\longrightarrow P'$ and completes the proof of Claim \ref{Ppushout}.

\

We may now use Claim \ref{Ppushout} to complete the description of the functor $\mathcal P_{1}$ and verify that $\eta$ is indeed a natural transformation. Let $\gamma$ be a path in $\ex C(f')_{T}$ joining $y$ to $y'$ so that $\phi_{T}\gamma$ is non constant and joins $x$ to $x'$. Then the following diagram commutes
\[\begin{tikzcd} &  \mathcal P(y) \ar{dl} \rar{\mathcal P\phi} \dar{\mathcal P(\gamma) }&\mathcal P(x) \dar{\mathcal P(\phi_{T}\gamma)}
\\\mathcal P_{0}(y)\ar{dr}{\mathcal P_{0}\gamma}& \mathcal P(y')\rar {\mathcal P\phi} \dar &\mathcal P(x')\dar
\\&\mathcal P_{0}(y')\rar{\eta} &\mathcal P_{1}(x')\end{tikzcd}\] 
As $\phi_{T}\gamma$ is non constant, $y$ must be equal to $\phi^{-1}(x)$. Consider $\mathcal P_{1}(x')$ as a constant functor from $\phi_{T}^{-1}(x)=y$ and use the outer loop in the above commutative diagram to apply  Claim \ref{Ppushout}. There therefore exists a unique map
\[\mathcal P_{1}(x)\longrightarrow \mathcal P_{1}(x')\]
so that the following diagram commutes
\begin{equation}\label{p1}\begin{tikzcd} &  \mathcal P(y) \ar{dl} \rar{\mathcal P\phi} &\mathcal P(x) \dar{\mathcal P(\phi_{T}\gamma)}\ar{dl}
\\\mathcal P_{0}(y)\ar{dr}{\mathcal P_{0}\gamma}\ar{r}{\eta}& \mathcal P_{1}(x)\ar[dashed]{dr}{\exists !} &\mathcal P(x')\dar
\\&\mathcal P_{0}(y')\rar{\eta} &\mathcal P_{1}(x')\end{tikzcd}\end{equation} 

 Define $\mathcal P_{1}(\phi_{T}\gamma)$ to be this map $\mathcal P_{1}(x)\longrightarrow \mathcal P_{1}(x')$. As every other map in the above diagram commutes with projection to $P$, Claim \ref{Ppushout} implies that $\mathcal P_{1}(\phi_{T}\gamma)$ also commutes with projection to $P$.
 \[\begin{tikzcd}\mathcal P_{1}(x)\ar{rr}{\mathcal P_{1}(\phi_{T}\gamma)}\ar{dr}&&\mathcal P_{1}(x')\ar{dl}
 \\ & P\end{tikzcd}\]
 
The fact that $\phi_{T}\gamma$  is non constant implies that $\mathcal P_{0}(y)=P=\mathcal P_{1}(x)$. If $\mathcal P_{1}(x')=P$ too, then the above commutative diagram implies that $\mathcal P_{1}(\phi_{T}\gamma)$ is the identity. Otherwise,  the map $\mathcal P_{1}(\phi_{T}\gamma):\mathcal P_{1}(x)\longrightarrow \mathcal P_{1}(x')$ is the inclusion of $P$ as the face of $\mathcal P_{1}(x')$ which contains the image of $\mathcal P(\phi_{T}\gamma)$ in $\mathcal P(x')\subset\mathcal P_{1}(x')$. For non constant paths $\gamma$ which are not in the image of $\phi_T$, we may similarly define $\mathcal P_{1}(\gamma)$ to be the inclusion of $\mathcal P_{1}(\gamma(0))=P$ as the face of $\mathcal P_{1}(\gamma(1))$ which contains the image of $\mathcal P(\gamma)$. This defines $\mathcal P_{1}$ as a functor.

 We have already observed that $\eta$ restricted to the inverse image of any object in $\ex C(f)_{T}$ is a natural transformation. The bottom left hand square of  diagram (\ref{p1}) then implies that $\eta$ is a natural transformation.

\smallskip

To describe $\mathcal P_{1}$ as an extension of the tropical structure of $f$ we must also describe a natural transformation $\mathcal P_{1}f$ from $\mathcal P_{1}$ to $\mathcal P\circ f_{T}$. Consider the commutative diagram of functors from $\ex C(f')_{T}$ to the category of integral-affine polytopes:
\[\begin{tikzcd}\mathcal P\rar \dar &\mathcal P\circ \phi_{T} \dar{\mathcal Pf \circ \phi_{T}}
\\ \mathcal P_{0} \rar{\mathcal P_{0}f'} &\mathcal P\circ f_{T}\circ \phi_{T}\end{tikzcd}\]

Claim \ref{Ppushout} then gives us for every object $x$ in $\ex C(f)_{T}$ a unique map 
\[\mathcal P_{1}f:\mathcal P_{1}(x)\longrightarrow \mathcal P(f(x))\]
so that  for any $y$ in $\phi_{T}^{-1}(x)$, the following diagram commutes
\begin{equation}\label{p1f}\begin{tikzcd}\mathcal P_{1}(x)\ar{dr}{\mathcal P_{1}f}&\mathcal P(x)\lar \dar{\mathcal Pf}
\\\mathcal P_{0}(y) \uar{\eta} \rar{\mathcal P_{0}f'}&\mathcal P(f(x))
 \end{tikzcd}\end{equation}
 We already know that every map apart from $\mathcal P_{1}f$ in the above diagram comes from a natural transformation.  Let $\gamma$ be a path from $y$ to $y'$ and $\phi_{T}(\gamma)$ be a path from $x$ to $x'$, and consider the diagram 
 \[\begin{tikzcd} \mathcal P_{1}(x)\ar{dr}{\mathcal P_{1}f}\ar[bend left]{rr}{\mathcal P_{1}(\phi_{T}(\gamma))}&\mathcal P(x)\lar \dar{\mathcal Pf} \ar[bend left]{rr}{\mathcal P(\phi_{T}(\gamma))} &  \mathcal P_{1}(x') \ar{dr}{\mathcal P_{1}f} &\mathcal P(x')\lar \dar{\mathcal Pf}
\\\mathcal P_{0}(y) \uar{\eta} \rar{\mathcal P_{0}f'}\ar[bend right]{rr}[swap]{\mathcal P_{0}(\gamma)} &\mathcal P(f(x))\ar[bend right]{rr}[swap]{\mathcal P(f_{T}(\phi_{T}(\gamma)))}
&\mathcal P_{0}(y') \uar{\eta} \rar{\mathcal P_{0}f'}&\mathcal P(f(x'))\end{tikzcd}\]
The commutativity of the squares corresponding to maps which we know are natural transformations implies that the square which should commute if $\mathcal P_{1}f$ is a natural transformation commutes when $\mathcal P_{1}(x)$ is restricted to the image of either $\mathcal P(x)$ or $\mathcal P_{0}(y)$. As these two subpolytopes span $\mathcal P_{1}(x)$, it follows that the square required for $\mathcal P_{1}f$ to be a natural transformation commutes. As any path in $\ex C(f)_{T}$ may be written as a composition of paths in the image of $\phi_{T}$ or their inverses, it follows that $\mathcal P_{1}f$ is a natural transformation.
 
The diagram (\ref{p1f}) above therefore corresponds to a commutative diagram of natural transformations. It follows that  $\mathcal P_{1}$ is an extension of the tropical structure of $f$, and  $\mathcal P_{1}f\circ \eta=\mathcal P_{0}f'$.

\smallskip

We have now constructed the required commutative diagram. 
\[\begin{tikzcd}[row sep=small](\ex C(f')_{T},\mathcal P_{0})\arrow[bend left=15]{drr}{(f'_{T},\mathcal P_{0}f')}\ar{d}{(\phi_{T},\eta)}
\\ (\ex C_{T},\mathcal P_{1})\ar{rr}{(f_{T},\mathcal P_{1}f)}&&(\ex B_{T},\mathcal P)
\end{tikzcd}\]

By the universal property of \(\mathcal P_u\), the extension \(\mathcal P_f\to\mathcal P_1\) determines a unique map of extended tropical structures \(\mathcal P_1\to\mathcal P_u\). We therefore obtain the following commutative diagram.
\[\begin{tikzcd}[row sep=small](\ex C(f')_{T},\mathcal P_{0})\arrow[bend left=15]{drr}{(f'_{T},\mathcal P_{0}f')}\ar{d}
\\ (\ex C_{T},\mathcal P_{1})\dar \ar{rr}{(f_{T},\mathcal P_{1}f)}&&(\ex B_{T},\mathcal P)
\\(\ex C_{T},\mathcal P_{u})\arrow[bend right=15]{rru}[swap]{(f_{T},\mathcal P_{u}f)} \end{tikzcd}\]
This gives the tropical part of our required map $\Phi$ from a neighborhood of $\ex C(f')$ in $\ex C(\hat h)$ to $\ex C(\hat f)$. We may now construct $\Phi$ as in the proof of Lemma \ref{uts map}.

\end{proof}

The following theorem is used in~\cite{evc} to construct a concrete local model for the moduli stack of stable curves in \(\hat{\ex B}\).
 
 \begin{thm}\label{G uts family}
For any stable curve \(f\) in \(\hat{\ex B}\) with domain not equal to \(\ex T\), there exists a family of curves \(\hat f\) containing \(f\), with universal tropical structure, and a group \(G\) of automorphisms of \(\hat f\), such that:
\begin{enumerate}
\item \(G\) acts freely and transitively on the set of maps of \(f\) into \(\hat f\);
\item there is a unique stratum \(\ex F_0\) of \(\ex F(\hat f)\) containing the image of all maps \(f\to\hat f\), and \(\totl{\ex F_0}\) is a single point;
\item the action of \(G\) on \(\totl{\ex C(\hat f)}\), restricted to the inverse image of \(\totl{\ex F_0}\), is effective, so \(G\) may be regarded as a subgroup of the automorphism group of \(\totl f\).
\end{enumerate}
\end{thm}
 
 \begin{proof}

Let \(\mathcal P_u\) be the universal extension of the tropical structure of \(f\), with base polytope \(P_u\). Lemma~\ref{uts family} constructs a family of curves \(\hat f_0\) containing \(f\), with universal tropical structure at \(f\). By restricting this family if necessary, we may assume that \(\ex F(\hat f_0)\) is an open subset of \(\et m{P_u}\), and that the stratum containing the image of \(f\) corresponds to the interior \(P_u^\circ\subset P_u\). Lemma~\ref{open uts} then implies that \(\hat f_0\) has universal tropical structure on this restricted family. In particular, any map of \(f\) into \(\hat f_0\) has image in the stratum \(\et m{P_u^\circ}\subset \et m{P_u}\).

The proof has four steps. First, we use Lemma~\ref{uts map} to construct a group \(G\) acting freely and transitively on the set of maps of \(f\) into \(\hat f_0\), restricted over \(\et m{P_u^\circ}\). Second, we show that \(G\) acts effectively on \(\totl f\), and hence may be regarded as a subgroup of the automorphism group of \(\totl f\). Third, we extend this \(G\)-action from the inverse image of \(\et m{P_u^\circ}\) to a neighbourhood in the total space of domains. Finally, we modify \(\hat f_0\), without changing it over \(\et m{P_u^\circ}\), so that the resulting map is \(G\)-invariant.

\smallskip
\noindent\emph{Constructing \(G\) over the interior stratum.}

Fix one inclusion of \(f\) into \(\hat f_0\). Given any other inclusion of \(f\) into \(\hat f_0\), Lemma~\ref{uts map}, applied to the two resulting families containing \(f\), gives a map
\[
        \Phi:\ex C(\hat f_0)\longrightarrow \ex C(\hat f_0)
\]
after restricting to a neighbourhood of \(f\), such that the distance between \(\hat f_0\) and \(\hat f_0\circ\Phi\) is bounded in any metric on \(\hat{\ex B}\). Over the inverse image of \(\et m{P_u^\circ}\), this bounded-distance condition implies
\[
        \hat f_0\circ\Phi=\hat f_0 .
\]
Moreover, the tropical part of \(\Phi\) is uniquely determined by the universal property of \(\mathcal P_u\). Since the smooth part of \(\et m{P_u^\circ}\) is a single point, this determines \(\Phi\) uniquely over the inverse image of \(\et m{P_u^\circ}\).

Thus any two inclusions of \(f\) into \(\hat f_0\) are exchanged by a unique automorphism of \(\hat f_0\) restricted to the inverse image of \(\et m{P_u^\circ}\). Let \(G\) be the group of these automorphisms. By construction, \(G\) acts freely and transitively on the set of maps of \(f\) into \(\hat f_0\).

\smallskip
\noindent\emph{Effectivity on the smooth part.}

The smooth part of the inverse image of \(\et m{P_u^\circ}\) in \(\ex C(\hat f_0)\) is equal to the smooth part of \(\ex C(f)\). Hence \(G\) acts by automorphisms of \(\totl f\). We claim that this action is effective.

Remark~\ref{pu char} implies that the interior stratum \(\et m{P_u^\circ}\) is detected by the tropical positions of a finite collection of marked point sections, together with the edge-length data recording how the smooth strata of the domain are glued at nodes. Therefore, two inclusions of \(f\) into \(\hat f_0\) with the same induced map on smooth parts have the same image in \(\et m{P_u^\circ}\): they have the same labelling of marked points and nodes, and the same corresponding tropical deformation data.

It follows that any element of \(G\) acting trivially on \(\totl{\ex C(f)}\) is an automorphism of the curve \(f\) itself. The only automorphism of \(f\) acting trivially on \(\totl{\ex C(f)}\) is the identity. Since \(G\) acts freely on the set of maps of \(f\) into \(\hat f_0\), the action of \(G\) on \(\totl{\ex C(f)}\) is effective. We may therefore regard \(G\) as a subgroup of the automorphism group of \(\totl f\). In particular, \(G\) is finite.

For curves in a smooth manifold, \(G\) is the automorphism group of \(\totl f\). In general, \(G\) is only a subgroup of the automorphism group of \(\totl f\).

\smallskip
\noindent\emph{Extending the action to the family of domains.}

We next extend the action of \(G\) to the total space \(\ex C(\hat f_0)\). In the construction of \(\ex C(\hat f_0)\) in the proof of Lemma~\ref{uts family}, we chose coordinate charts \(U_i\) on \(\ex C(f)\), and then extended them to coordinate charts \(\tilde U_i\) on \(\ex C(\hat f_0)\). We may assume that the action of \(G\) preserves the smooth parts of these charts, so that the smooth part of each \(U_i\) is sent to the smooth part of some \(U_j\). For instance, one can achieve this by choosing a \(G\)-invariant metric on \(\totl{\ex C(f)}\) in the conformal class determined by the complex structure, and then choosing the charts of type~\ref{t1} and type~\ref{t2} to be the points within a fixed distance of the relevant node or marked point.

The smooth part of \(\tilde U_i\), restricted to the inverse image of \(\et m{P_u^\circ}\), is the smooth part of \(U_i\). Thus, over \(\et m{P_u^\circ}\), every element of \(G\) sends \(\tilde U_i\) isomorphically into some \(\tilde U_j\). These maps are equivariant in the sense of~\cite{cem}: they send smooth monomials to smooth monomials. By Remark~\ref{equivariant coordinates}, the coordinate system \(\{\tilde U_I\}\) is equivariant, so each such map over \(\et m{P_u^\circ}\) has a unique extension to an equivariant map between the corresponding coordinate charts.

Concretely, if \(U_i\) is a chart of type~\ref{t1} or type~\ref{t2}, then the extended map
\[
        \Phi:\tilde U_i\longrightarrow \tilde U_j
\]
is the unique monomial map restricting to the given map over the inverse image of \(\et m{P_u^\circ}\). If \(U_i\) is a chart of type~\ref{t3}, then \(\tilde U_i=U_i\times\ex F(\hat f_0)\), and the extended map is uniquely the product of a map \(U_i\to U_j\) with a monomial automorphism of \(\ex F(\hat f_0)\), again restricting to the given map over the inverse image of \(\et m{P_u^\circ}\). The definition of the fiberwise complex structure on the charts \(\tilde U_i\) makes these equivariant maps automatically fiberwise complex.

The uniqueness of these equivariant extensions implies that the automorphisms in \(G\), initially defined over the inverse image of \(\et m{P_u^\circ}\), extend to an action of \(G\) by automorphisms of \(\ex C(\hat f_0)\), after restricting \(\ex C(\hat f_0)\) to a neighbourhood of that inverse image if necessary.

\smallskip
\noindent\emph{Averaging the map to make it \(G\)-invariant.}

It remains to modify the map \(\hat f_0\), away from the inverse image of \(\et m{P_u^\circ}\), so that it is preserved by the \(G\)-action. Let \(U\) be a \(G\)-invariant open subset of \(\ex C(\hat f_0)\) such that each connected component of \(U\) is simply connected and intersects the inverse image of \(\et m{P_u^\circ}\) in a connected subset. Assume also that \(\hat f_0(U)\) is contained in a coordinate chart
\[
        V=\mathbb R^n\times \et kP
\]
on \(\hat{\ex B}\), such that \(T_{\mathrm{vert}}\hat{\ex B}\) over \(V\) is spanned by a fixed subset of the standard vector fields
\[
        \left\{\frac{\partial}{\partial x_i},\ \tilde z_j\frac{\partial}{\partial z_j}\right\}.
\]

For \(g\in G\), the map \(\hat f_0\circ g\) on \(U\) is obtained from \(\hat f_0\) by exponentiating a section \(v_g\) of \(\hat f_0^*T_{\mathrm{vert}}\hat{\ex B}\), using the locally defined connection on \(V\) which preserves the standard basis vector fields. The section \(v_g\) is uniquely determined by requiring that it vanish on the inverse image of \(\et m{P_u^\circ}\). Since \(G\) is finite, we may average these sections over \(G\). Replacing \(\hat f_0\) on \(U\) by the map obtained by exponentiating the average of the \(v_g\) gives a \(G\)-invariant map \(U\to V\), and this new map agrees with \(\hat f_0\) on the inverse image of \(\et m{P_u^\circ}\).

Using a cutoff function, we interpolate between \(\hat f_0\) and this averaged map. This gives a new map which is \(G\)-invariant on a compactly contained open subset \(U_0\subset U\), and which agrees with \(\hat f_0\) on the inverse image of \(\et m{P_u^\circ}\). Moreover, if \(\hat f_0\) was already \(G\)-invariant on some connected component of \(U\) intersecting the inverse image of \(\et m{P_u^\circ}\), then all the sections \(v_g\) vanish there, so the modification does not change \(\hat f_0\) on that component.

Repeating this averaging construction finitely many times, for a \(G\)-invariant collection of such open sets, we obtain a map \(\hat f\) which is \(G\)-invariant on a \(G\)-invariant neighbourhood of the inverse image of \(\et m{P_u^\circ}\). Restrict the family to this neighbourhood.

The resulting family \(\hat f\) agrees with \(\hat f_0\) over the inverse image of \(\et m{P_u^\circ}\). Hence \(G\) still acts freely and transitively on the set of maps \(f\to\hat f\), and its action on the smooth part over the corresponding stratum is effective. The only stratum of \(\ex F(\hat f)\) containing the image of a map \(f\to\hat f\) is the stratum \(\ex F_0\) corresponding to \(P_u^\circ\), and \(\totl{\ex F_0}\) is a single point. Finally, the modification from \(\hat f_0\) to \(\hat f\) is by exponentiating vertical vector fields and interpolating, so it does not change the relevant tropical structure. Since \(\hat f_0\) has universal tropical structure, the restricted \(G\)-invariant family \(\hat f\) also has universal tropical structure. This proves the theorem.

\end{proof}

\bibliographystyle{plain}
\bibliography{../ref}
 \end{document}